\documentclass[fleqn,a4paper,10pt]{amsart}

\usepackage{graphicx}
\usepackage{tikz}
\usetikzlibrary{shapes,decorations.pathmorphing}
\usepackage{amsmath}
\usepackage{amssymb,amsthm}
\usepackage{comment}
\usepackage{hyperref}

\theoremstyle{definition}

\title[Statistics of hull perimeters]{Some results on the statistics of hull perimeters in large planar triangulations and quadrangulations}
\author{Emmanuel Guitter}
\address{Institut de Physique Th\'eorique, CEA, IPhT, 91191 Gif-sur-Yvette, France, CNRS, UMR 3681}
\email{emmanuel.guitter@cea.fr}

\begin{document}
\maketitle

\begin{abstract}
The hull perimeter at distance $d$ in a planar map with two marked vertices at distance $k$ from each other is the length of the closed curve separating these two vertices and lying at distance $d$ from the first one ($d<k$). We study the statistics of hull perimeters in large random planar triangulations and quadrangulations as a function of both $k$ and $d$. Explicit expressions for the probability density
of the hull perimeter at distance $d$, as well as for the joint probability density of hull perimeters at distances $d_1$ and $d_2$, are obtained
in the limit of infinitely large $k$. We also consider the situation where the distance $d$ at which the hull perimeter is measured corresponds 
to a finite fraction of $k$. The various laws that we obtain are identical for triangulations and for quadrangulations, up to a global rescaling.
Our approach is based on recursion relations recently introduced by the author which determine the 
generating functions of so-called slices, i.e.\ pieces of maps with appropriate distance constraints. It is indeed shown that the map
decompositions underlying these recursion relations are intimately linked to the notion of hull perimeters and provide a simple way to 
fully control them.
\end{abstract}

\section{Introduction}
\label{sec:introduction}
\begin{figure}
\begin{center}
\includegraphics[width=10cm]{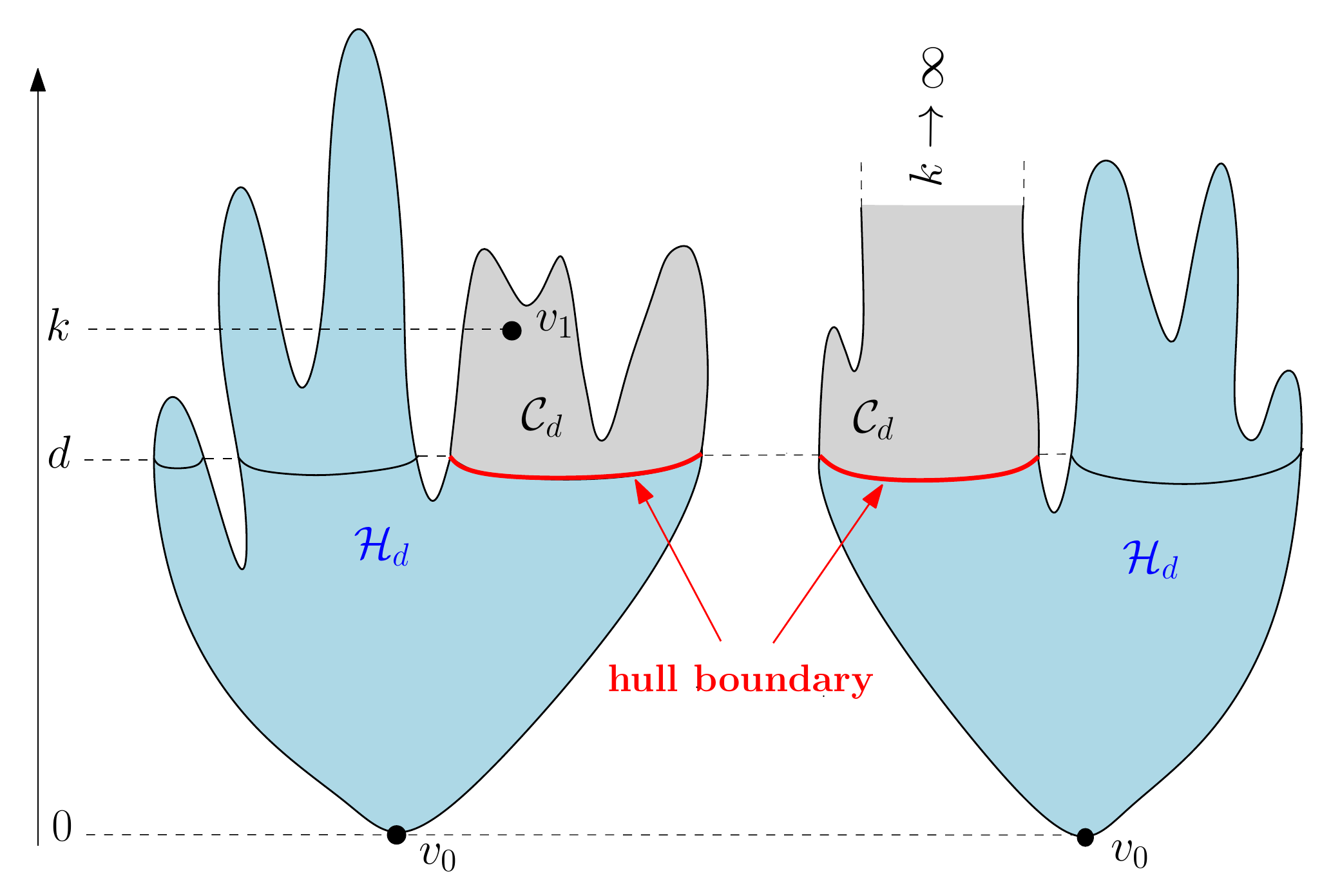}
\end{center}
\caption{Left: A schematic picture of the hull boundary in a planar map with an origin $v_0$ and a second marked vertex $v_1$ at distance $k$
from $v_0$. The map is represented so that vertices appear at a height equal to their distance from $v_0$ . The vertices
at distance $d$ from $v_0$ form a number of closed curves at height $d$, one of which separates $v_0$ from $v_1$ and constitutes
the hull boundary.  Right: When the map is infinite, exactly one of the connected components formed by vertices at distance larger that $d$ 
is infinite. Sending $k\to \infty$ ensures that $v_1$ belongs to this component. This heuristic view will be made rigorous by precise definitions
along the paper.}
\label{fig:hullboundary}
\end{figure}
Understanding the statistics of random planar maps, i.e.\ connected graphs embedded on the sphere, as well as their various continuous limits,
such as the Brownian map \cite{Miermont2013,legall2013} or the Brownian plane \cite{CLG13}, is a very active field of both 
combinatorics and probability theory.
In this paper, we study the \emph{statistics of hull perimeters} in large planar maps, a problem which may heuristically be understood as follows: consider a planar map $\mathcal{M}$ of some type (in the following, we shall restrict our analysis to the case of triangulations and quadrangulations) with two marked vertices, an origin vertex $v_0$ and a second distinguished vertex $v_1$ at graph distance $k$ from $v_0$, for some $k \geq 2$. 
Consider now, for some $d$ strictly between $0$ and $k$, the \emph{ball} of radius $d$ which, so to say, is the part of the map at graph distance
less than $d$ from the origin\footnote{Several prescriptions may be used to precisely define the ball, each leading to a slightly different
notion of hull.}. This ball has a boundary made in general of several closed lines, each line linking vertices at distance of order $d$ from $v_0$ and separating $v_0$ from a connected domain where all vertices are at distance larger than $d$ (see figure \ref{fig:hullboundary}). One of these domains 
$\mathcal{C}_d$ contains the second distinguished vertex $v_1$ and we may define the \emph{hull} of radius $d$ as the domain 
$\mathcal{H}_d=\mathcal{M}\setminus\mathcal{C}_d$,  namely the union of the ball of radius $d$ itself and of all the connected 
domains at distance larger than $d$ which do not contain $v_1$ (see figure \ref{fig:hullboundary} where $\mathcal{H}_d$
is represented in light blue). The hull boundary is then the boundary of $\mathcal{H}_d$ (which is also
that of $\mathcal{C}_d$), forming a closed line at distance $d$ from $v_0$ and separating $v_0$ from $v_1$. 
The length of this boundary is called the \emph{hull perimeter at distance $d$} and will be denoted by $\mathcal{L}(d)$ in the following.

The purpose of this paper is to study, as a function of both $k$ and $d$, the statistics of the hull perimeter $\mathcal{L}(d)$ within uniformly
drawn planar maps of some given type (here triangulations and quadrangulations) equipped with a randomly chosen pair of marked vertices at distance $k$ from one another. Even though the combinatorics developed in this paper allows us to keep the size $N$ ($=$ number of faces) of the
maps finite, explicit statistical laws will be presented only in the limit of \emph{infinitely large maps}, namely when $N\to \infty$. In this case, 
it is  expected that exactly one of the components outside the ball of radius $d$ is infinite. In particular,  
sending $k \to \infty$ allows us to enforce that $v_1$ belongs to this infinite component so that the hull of radius $d$
no longer depends on $v_1$ in this case. This limits describes a slightly simpler notion of hull boundary for vertex-pointed 
infinite planar maps, namely the line at distance $d$ from the origin $v_0$ separating this origin from infinity 
(see figure \ref{fig:hullboundary}-right). 

The question of the hull perimeter statistics was already addressed in several papers \cite{Krikun03,Krikun05,CLG14a,CLG14b}.
For a given choice of map ensemble, the above heuristic presentation may be transformed into a well-defined statistical problem by a 
rigorous definition of the hull boundary at distance $d$ within the maps at hand. 
Several prescriptions may be adopted and our precise choice will be detailed in Section \ref{sec:definitions}. This choice is different
from the prescriptions used in \cite{Krikun03,Krikun05,CLG14a,CLG14b} and is thus expected to give different results for the hull perimeter statistics at finite $d$ and $k$. 
In the limit of large $d$ and $k$ however, all prescriptions should eventually yield the same universal laws (up to some possible finite rescalings). 
This assumption is corroborated by our results on the hull perimeter probability density which precisely reproduce the expressions of
\cite{Krikun03,Krikun05,CLG14a,CLG14b}, as displayed in eqs.~\eqref{eq:CLG} and \eqref{eq:Pinf} below. The laws that we find for large $d$ and $k$ are the same for triangulations and for quadrangulations, up to a global scale change, and we expect that they should emerge for other families 
of maps as well. 
\vskip .3cm
The paper is organized as follows: we start by giving in Section~\ref{sec:definitions} our precise definition of the hull boundary,
which we view as a particular dividing line drawn on some canonical ``slice" representation of the map at hand. The
construction of this line is slightly different for triangulations and for quadrangulations and mimics that discussed
by the author in \cite{G15a,G15b} in a related context. We then present our main results in Section~\ref{sec:mainresults},
namely explicit expressions for the probability density of the hull perimeter in (i) the regime of a large but finite $d$ and $k\to \infty$
and (ii) the regime of large $d$ and $k$ with a fixed ratio $d/k<1$, as well as for the joint probability density of the hull perimeters at
two large but finite distances $d_1$ and $d_2$ for $k\to\infty$. Sections~\ref{sec:recursion} and \ref{sec:largemaps} are devoted to the derivation
of our main results. We first recall in Section~\ref{sec:recursion} the existence of a recursion relation for the generating function of
the slices representing our maps and show how the map decomposition underlying this recursion may be related to our notion of hull boundary.
This allows us to obtain easily a number of explicit expressions for map generating functions with a control on hull perimeters.
The case of quadrangulations is discussed in Section~\ref{sec:quad} and that of triangulations in Section~\ref{sec:triang}.
Expansions of these generating functions are presented in Appendix B, which display the \emph{numbers} of quadrangulations
and triangulations with fixed $N$ ($=$ number of faces), $k$, $d$ and $\mathcal{L}(d)$ for the first allowed values.
We then extract in Section~\ref{sec:largemaps}  from the singularities of the generating functions
the desired hull perimeter probability densities for quadrangulations and triangulations with an infinitely large number $N$ of faces.  
The details of the technique are discussed in Section~\ref{sec:singularity} and we show in Section~\ref{sec:shortcut}
how to slightly simplify the calculations for large $k$. Some involved intermediate formulas are given in Appendix C.
The case of large $d$ and $k$ of the same order is
discussed in Section~\ref{sec:finitefraction}. All over the paper, explicit expressions are first obtained for the \emph{Laplace transforms}
of the various probability densities at hand. The final step consisting in taking the desired inverse Laplace transforms 
is discussed in Appendix A.  
We gather our concluding remarks in Section~\ref{sec:conclusion}.

\section{Summary of the results}
\label{sec:summary}
\subsection{Definition of the hull boundary in pointed-rooted triangulations and quadrangulations}
\label{sec:definitions}
A natural way to define the hull boundary is to first construct the ball of radius $d$: this requires deciding which \emph{faces} of the map are retained in the ball\footnote{As opposed to vertices which are simply characterized by their graph
distance from $v_0$, faces may be of different types
according to the graph distance of their incident vertices.} and many inequivalent choices may be adopted, each leading to a slightly different definition of hull. In \cite{Krikun03},  Krikun gave a particularly elegant prescription in the case of triangulations, later used 
in \cite{CLG14a}, which allowed him to relate the hull boundary statistics
to that of some ``time-reversed" branching process. An alternative way to define the ball and hull of radius $d$,
described in \cite{CLG14b}, is to use
the graph distance on the dual map as a way one to assign distances directly to the faces of the original
map. Here we shall use yet another prescription and  \emph{construct directly} the hull boundary at distance $d$ without having recourse to a preliminary construction of the ball of radius $d$. Our approach applies to both triangulations and quadrangulations and is based on a technique developed recently in \cite{G15a,G15b} to compute the twopoint function of these maps.

\begin{figure}
\begin{center}
\includegraphics[width=9cm]{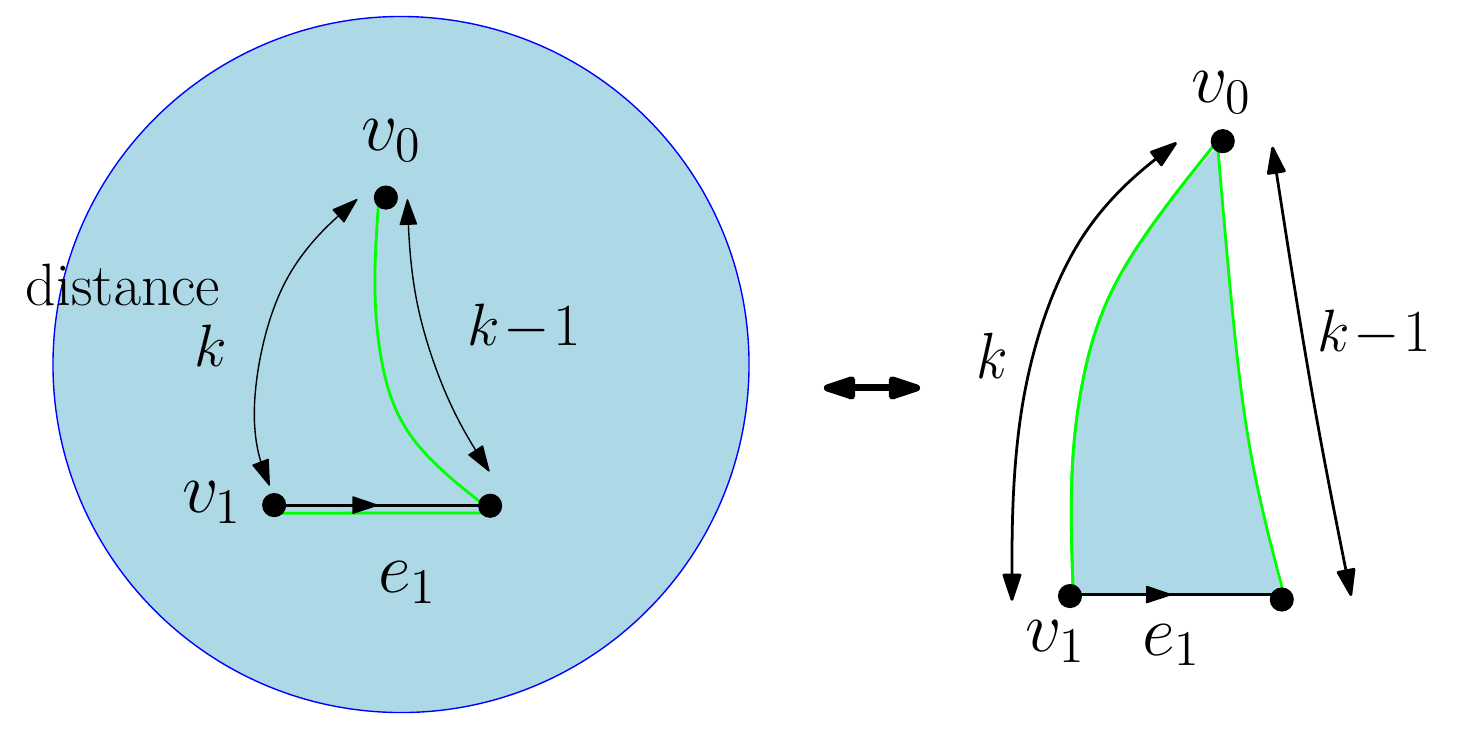}
\end{center}
\caption{A schematic picture of a $k$-pointed-rooted planar map (left), i.e.\ a map with a marked origin $v_0$ and a marked edge $e_1$ 
with extremities at distance $k$ and $k-1$ from $v_0$. We have drawn in green the leftmost shortest path (starting with $e_1$) from
$v_1$ to $v_0$. Cutting along this line and unwrapping the map creates a $k$-slice with apex $v_0$ and base $e_1$ (right)
characterized by the properties (1)-(4) of the text. The correspondence between $k$-pointed-rooted maps and $k$-slices
is a bijection.}
\label{fig:twopoint}
\end{figure}
Our starting point is an arbitrary \emph{$k$-pointed-rooted} planar triangulation (respectively quadrangulation) i.e.\ a planar map whose all faces have degree $3$ (respectively $4$)
endowed with a marked vertex $v_0$ (the origin vertex) and a marked edge $e_1$ (the root edge) oriented from a vertex $v_1$ at distance $k$
from $v_0$ towards a vertex at distance $k-1$ from $v_0$ (see figure \ref{fig:twopoint}-left), for some $k\geq 1$ (respectively $k\geq 2$). The 
use of a pointed-rooted map rather than a simple vertex bi-pointed map (with two marked vertices $v_0$ and $v_1$) is a standard procedure which 
highly simplifies the underlying combinatorics. Note that, by definition, not all edges leaving a vertex $v_1$ at distance $k$ from $v_0$
can serve as root edge for our pointed-rooted map (since we impose that $e_1$ necessarily points towards a vertex at distance $k-1$ from $v_0$)
but that, for each $v_1$ at distance $k$, at least one such edge exists.

It is well-known that 
any $k$-pointed-rooted planar map, as defined above, may be bijectively transformed into a so-called \emph{$k$-slice} by cutting the map along the \emph{leftmost}
shortest path\footnote{This path is the sequence of edges obtained by taking $e_1$ as first step and then, at each encountered vertex at distance $\ell$ from $v_0$, picking the
leftmost edge leading from this vertex to a vertex at distance $\ell-1$ from $v_0$, until the path eventually reaches $v_0$.} from $v_1$ to $v_0$ and then unwrapping the map (see figure \ref{fig:twopoint}-right). A $k$-slice is
a planar map whose all faces have a fixed degree $\delta$ (with $\delta=3$ if the map before cutting was a triangulation or $\delta=4$ if it was a quadrangulation), except the root face (which is the face on the right of $e_1$ after cutting
and unwrapping) which has degree $2k$.  A $k$-slice is characterized by the following properties:
\begin{enumerate}
 \item the left boundary of a $k$-slice, which is formed by the $k$ edges incident to its root face lying between $v_1$ and $v_0$
 clockwise around the rest of the map is a shortest path between $v_1$ and $v_0$ within the $k$-slice;
 \item the  right boundary of a $k$-slice, which is formed by the $k-1$ edges incident to its root face lying between the endpoint of $e_1$ and $v_0$
 counterclockwise around the rest of the map is a shortest path between these two vertices within the $k$-slice; 
 \item the right boundary is the \emph{unique} shortest path between its endpoints within the $k$-slice;
 \item the left and right boundaries do not meet before reaching $v_0$.
 \end{enumerate}
The vertex $v_0$ is called the \emph{apex} and the edge $e_1$ the \emph{base} of the $k$-slice. Clearly, due to the particular choice of cutting line, the cutting procedure 
applied to a $k$-pointed-rooted map creates a $k$-slice. Given this $k$-slice, 
the original $k$-pointed-rooted map is reconstructed by gluing the left boundary of the $k$-slice to the union of its base $e_1$ and its right boundary 
in the unique way which preserves the distances to $v_0$. 
The transformation between $k$-pointed-rooted maps and $k$-slices is a bijection.  
More simply, the $k$-slice with apex $v_0$ and base $e_1$ may be viewed as a canonical representation of the associated $k$-pointed-rooted map
with origin $v_0$ and root edge $e_1$.
 
\subsubsection{Definition of the hull boundary in a pointed-rooted quadrangulation.}
\begin{figure}
\begin{center}
\includegraphics[width=11cm]{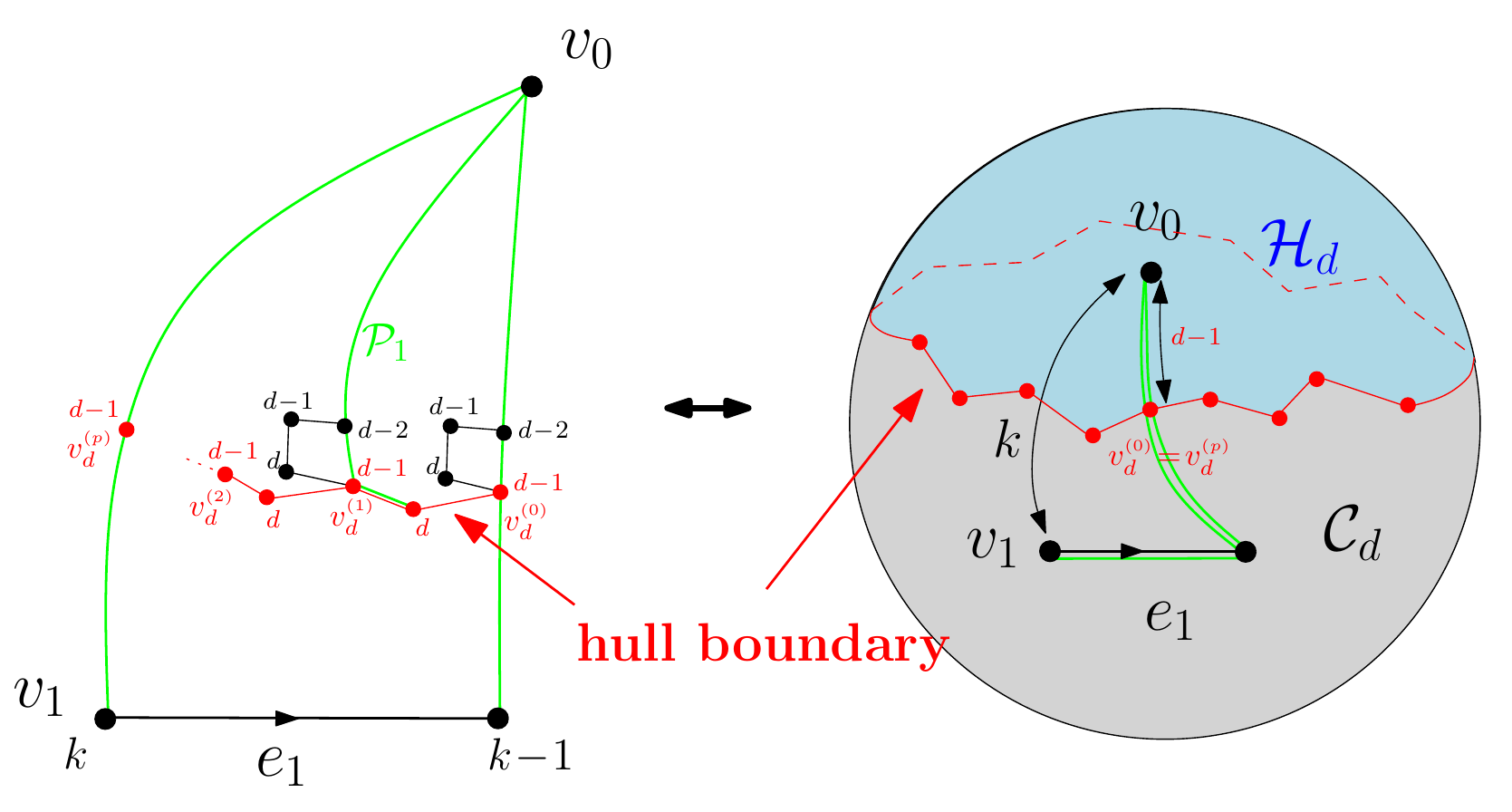}
\end{center}
\caption{The construction of the hull boundary for quadrangulations (see text for details). The hull boundary is first constructed as an open line 
on the associated $k$-slice, then transformed into a  simple closed curve on the quadrangulation by re-gluing.}
\label{fig:dividing}
\end{figure}
Let us now define the hull boundary of our $k$-pointed-rooted map. We start with the case of quadrangulations and perform
our construction on the associated $k$-slice. Our construction is in all point similar to that discussed 
in \cite{G15b}. Given $k\geq 3$ and $d$ in the range $2\leq d\leq k-1$, we start from the (unique) vertex $v_d^{(0)}$ of the right boundary of the 
$k$-slice at distance $d-1$ from $v_0$ and consider the face on the left of the edge connecting $v_d^{(0)}$ to its neighbor at distance $d-2$ from $v_0$ along the right boundary. The last two vertices incident to the face are necessarily distinct from the first two and are at distance $d$ and $d-1$ 
respectively from $v_0$ (see figure \ref{fig:dividing}). 
This is a direct consequence of the property (3) of the right boundary (see \cite{G15b} for a detailed argument). The vertex $v_d^{(0)}$ is therefore incident to at least one edge leading to a neighboring vertex at distance $d$ from $v_0$ itself incident to one edge leading to a vertex at distance $d-1$ from $v_0$ and distinct from $v_d^{(0)}$. These two edges define a \emph{2-step path} of ``type" $d-1\to d\to d-1$ with distinct endpoints. Pick
the \emph{leftmost} such 2-step path from $v_d^{(0)}$ 
and call $(e^{(0)}_d,f^{(0)}_d$) the corresponding pair of successive edges, leading to a vertex $v^{(1)}_d$
at distance $d-1$ from $v_0$ and different from $v_d^{(0)}$. We can then draw the leftmost shortest path $\mathcal{P}_1$ from $v_d^{(1)}$ to $v_0$ and consider the face 
on the left of the edge connecting $v_d^{(1)}$ to its neighbor on $\mathcal{P}_1$ at distance $d-2$ from $v_0$. Repeating the argument allows us to construct
a leftmost 2-step path $(e^{(1)}_d,f^{(1)}_d)$ from $v_d^{(1)}$ to yet another distinct vertex $v_d^{(2)}$ and so on. As explained in \cite{G15b}, the line formed by the successive 2-step paths
$(e^{(i)}_d,f^{(i)}_d)$, $i=0,1,\cdots$ cannot form loops in the $k$-slice and necessarily ends after $p$ iterations at the (unique) vertex $v_d^{(p)}$ at
distance $d-1$ from $v_0$ lying on the left boundary of the $k$-slice\footnote{Note that the vertex preceding $v_d^{(p)}$ on the line,
at distance $d$ from $v_0$ may lie either strictly inside the $k$-slice or on its left boundary.} (see \cite{G15b} for a detailed argument).
This line forms our hull boundary at distance $d$ from $v_0$. Indeed, upon re-gluing the  $k$-slice into a $k$-pointed-rooted quadrangulation, 
we identify $v_d^{(0)}$ and $v_d^{(p)}$ and the line forms a simple closed curve visiting alternatively vertices at distance $d-1$ and $d$ from $v_0$ and separating $v_0$ from 
$v_1$ (see figure \ref{fig:dividing}). Clearly, all the vertices in the domain lying on the same side of the line as $v_1$ are at a distance larger than or equal 
to $d-1$ (and which can be equal to $d-1$ on the line only) and this domain constitutes  what we called $\mathcal{C}_d$ in the introduction.
As for the domain lying on the same side of the line as $v_0$, it contains all the vertices at distance less than or
equal to $d-1$ from $v_0$ (together with other vertices at arbitrary distance) and constitutes the hull $\mathcal{H}_d$.

The hull perimeter $\mathcal{L}(d)$ is the length $2p$ of the line above. 
For convenience, we decide to extend our definition of the hull boundary to the case $d=1$ and $k\geq 2$ by taking the convention that it is
then reduced to the simple vertex $v_0$ and has length $0$ accordingly, i.e.:
\begin{equation*}
\mathcal{L}(1)=0\quad \hbox{for quadrangulations} \ .
\end{equation*}

\subsubsection{Definition of the hull boundary in a pointed-rooted triangulation.}
\begin{figure}
\begin{center}
\includegraphics[width=11cm]{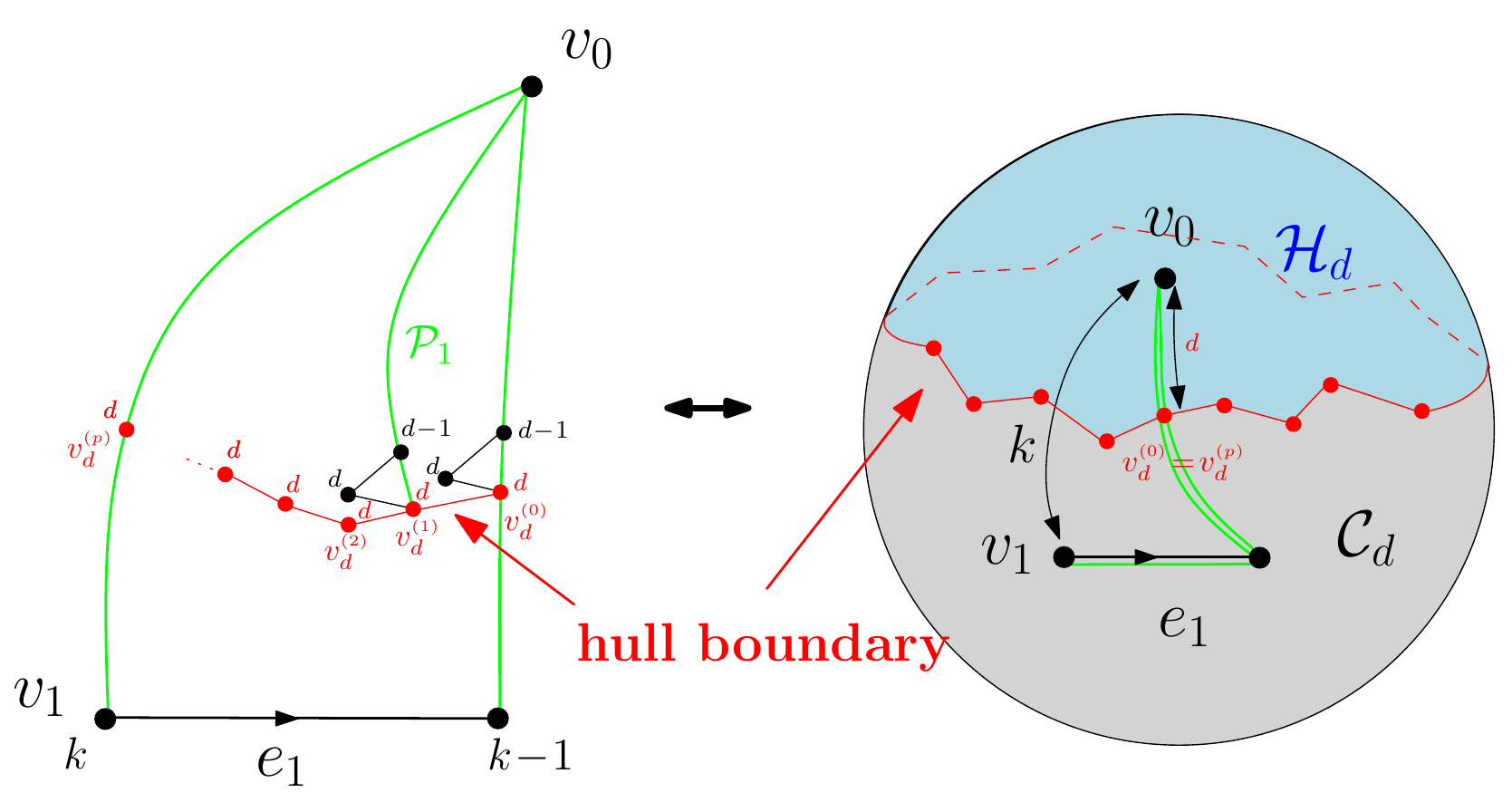}
\end{center}
\caption{The construction of the hull boundary for triangulations (see text for details). The hull boundary is first constructed as an open line 
on the associated $k$-slice, then transformed into a  simple closed curve on the triangulation by re-gluing.}
\label{fig:dividing2}
\end{figure}
The case of triangulations is slightly simpler and the construction follows that of
\cite{G15a}. We start from the $k$-slice associated with a given $k$-pointed-rooted triangulation, with now $k\geq 2$. 
Given $d$ in the range $1\leq d \leq k-1$, we look at the (unique) vertex $v_d^{(0)}$ of the right boundary of the $k$-slice at distance $d$ from $v_0$
and consider the triangle on the left of the edge connecting $v_d^{(0)}$ to its neighbor at distance $d-1$ from $v_0$ along the right boundary.
The third vertex incident to the triangle is necessarily distinct from the first two and at distance $d$ from $v_0$ (see figure \ref{fig:dividing2}). 
This is again a direct
consequence of the property (3) of the right boundary (see \cite{G15a} for a detailed argument). The vertex $v_d^{(0)}$ is therefore incident
to at least one edge leading to a distinct neighbor at distance $d$ from $v_0$. Call $e^{(0)}_d$ the leftmost such edge and $v^{(1)}_d$ its
extremity different from $v_d^{(0)}$. We can then draw the leftmost shortest path $\mathcal{P}_1$ from $v_d^{(1)}$ to $v_0$ and consider the triangle 
on the left of the edge connecting $v_d^{(1)}$ to its neighbor on $\mathcal{P}_1$ at distance $d-1$ from $v_0$. 
Repeating the argument allows us to construct
a leftmost edge $e^{(1)}_d$ to yet another distinct vertex $v_d^{(2)}$ and so on. As explained in \cite{G15a}, the line formed by the successive edges 
$e^{(i)}_d$, $i=0,1,\cdots$ cannot form loops in the $k$-slice and necessarily ends after $p$ iterations at the (unique) vertex $v_d^{(p)}$ at
distance $d$ from $v_0$ lying on the left boundary of the $k$-slice (see \cite{G15a} for a detailed argument).
This line forms our hull boundary at distance $d$ from $v_0$. Upon re-gluing the $k$-slice into a $k$-pointed-rooted triangulation, 
the line indeed forms a simple closed curve visiting only vertices at distance $d$ from $v_0$ and separating $v_0$ from 
$v_1$. Clearly, all the vertices in the domain lying on the same side of the line as $v_1$ are at a distance larger than or equal 
to $d$ (and which can be equal to $d$ on the line only): this domain constitutes  what we called $\mathcal{C}_d$ in the introduction.
The complementary domain $\mathcal{H}_d$, lying on the same side of the line as $v_0$, constitutes the hull and contains all the vertices at distance less than or
equal to $d$ from $v_0$ together with other vertices at arbitrary distance.

The hull perimeter $\mathcal{L}(d)$ is the length $p$ of the line above. 
Again, for convenience, we decide to extend our definition of the hull boundary to the case $d=0$ and $k\geq 1$ by taking the convention that it is
then reduced to the simple vertex $v_0$ and has length $0$, namely:
\begin{equation*}
\mathcal{L}(0)=0\quad \hbox{for triangulations} \ .
\end{equation*}

\subsection{Main results on the statistics of hull perimeters}
\label{sec:mainresults}
Having defined the hull perimeter $\mathcal{L}(d)$, our results concern the statistics of this perimeter in the ensemble of uniformly drawn 
$k$-pointed-rooted
quadrangulations (respectively triangulations) having a \emph{fixed number of faces} $N$, and for a \emph{fixed value of the parameter $k$}
(recall that our definition of 
$k$-pointed-rooted maps imposes not only that the first extremity $v_1$ of their root edge $e_1$ is at distance $k$ from the origin $v_0$
but also that the second extremity of $e_1$ is at distance $k-1$ from $v_0$). 
More precisely, we shall give explicit expressions \emph{in the limit $N\to \infty$ of this ensemble}. Note that, when sending $N\to \infty$,  
$k$ is kept finite (at least at a first stage) and \emph{does not scale with $N$}. 
This limit is called the \emph{local limit} of infinitely large quadrangulations (respectively triangulations).
We shall denote by $P_k(\{\cdot\})$ the probability of some event $\{\cdot\}$ and $E_k(\{\cdot\})$ the expectation value of some quantity 
$\{\cdot\}$ in this limit. 
When $d$ and $k$ themselves become large (recall that $d\leq k-1$), we find that the perimeter $\mathcal{L}(d)$ typically 
scales like $d^2$, so we are  naturally led to define the rescaled quantity:
\begin{equation*}
L(d)\equiv\frac{\mathcal{L}(d)}{d^2}\ .
\end{equation*}
Let us now present the main three results of this paper on the statistics of $L(d)$.
 
\subsubsection{Probability density for $L(d)$ when $k\to \infty$}
Our first result concerns the $k\to \infty$ limit, with probabilities and expectation values denoted by  $P_\infty(\{\cdot\})$
and $E_\infty(\{\cdot\})$. We insist here on the fact that,
although both $N$ and $k$ are sent to infinity, $k$ does not scale with $N$: in other words, we first send $N\to \infty$, and only then
send $k\to \infty$. As mentioned in the introduction, the hull perimeter $\mathcal{L}(d)$ may then be viewed as the length of the line ``at distance 
$d$" from $v_0$ separating $v_0$ from infinity. We find:
\begin{equation}
\lim_{d\to \infty} E_\infty (e^{-\tau L(d)})= \frac{1}{(1+c\, \tau)^{3/2}}\ ,
\label{eq:CLG}
\end{equation}
or equivalently (via a simple inverse Laplace transform):
\begin{equation}
\lim_{d\to \infty} P_\infty(L\leq L(d)< L+dL)= \frac{2}{\sqrt{\pi}}\frac{\sqrt{L}}{c^{3/2}}e^{-\frac{L}{c}}dL
\label{eq:Pinf}
\end{equation}
with $c$ taking a different value for quadrangulations and triangulations, namely:
\begin{equation}
\left\{\begin{matrix}
& c=\displaystyle{\frac{1}{3}}& \ \hbox{for quadrangulations}\ , \\
& & \\
& c=\displaystyle{\frac{1}{2}}&\ \hbox{for triangulations . \ \ \ \ }\\
\end{matrix}
\right.
\label{eq:valc}
\end{equation}
We recover here the precise form of the hull perimeter probability density found by Krikun \cite{Krikun03,Krikun05} and by   
Curien and Le Gall \cite{CLG14a,CLG14b}. The value of $c=1/3$ that we find for quadrangulations matches that of \cite{CLG14b}\footnote{The
correspondence with \cite{CLG14b} is $c=p/(4 h^2)$ where $p=2^{2/3}/3$ and $h=1/2^{2/3}$ for quadrangulations, leading to $c=1/3$
and $p=1/3^{1/3}$ and $h=1/(2\cdot 3^{1/6})$ for triangulations, leading to $c=1$. Although $p$ and $h$ are different, the same value $c=1$ is also obtained in \cite{CLG14b} for so-called triangulations of type II, which have no loops. These are the triangulations considered in
\cite{Krikun03} by Krikun, who also finds $c=1$, while he gets $c=1/2$ for quadrangulations \cite{Krikun05}.}
but is only $2/3$ of that found in \cite{Krikun05}. This suggest that our prescription and that of \cite{CLG14b} yield hull boundaries whose lengths 
are essentially the same, while the prescription used in \cite{Krikun05} creates hull boundaries which are larger by a factor $3/2$. 
Our value $c=1/2$ for triangulations is half that found in \cite{Krikun03,CLG14b}, suggesting that our prescription yields hull boundaries whose lengths are half those of the previous studies.   

Note that, we find in particular, 
\begin{equation*}
\lim_{d\to \infty} E_\infty (L(d))=\frac{3c}{2}\ .
\end{equation*}
\begin{figure}
\begin{center}
\includegraphics[width=9cm]{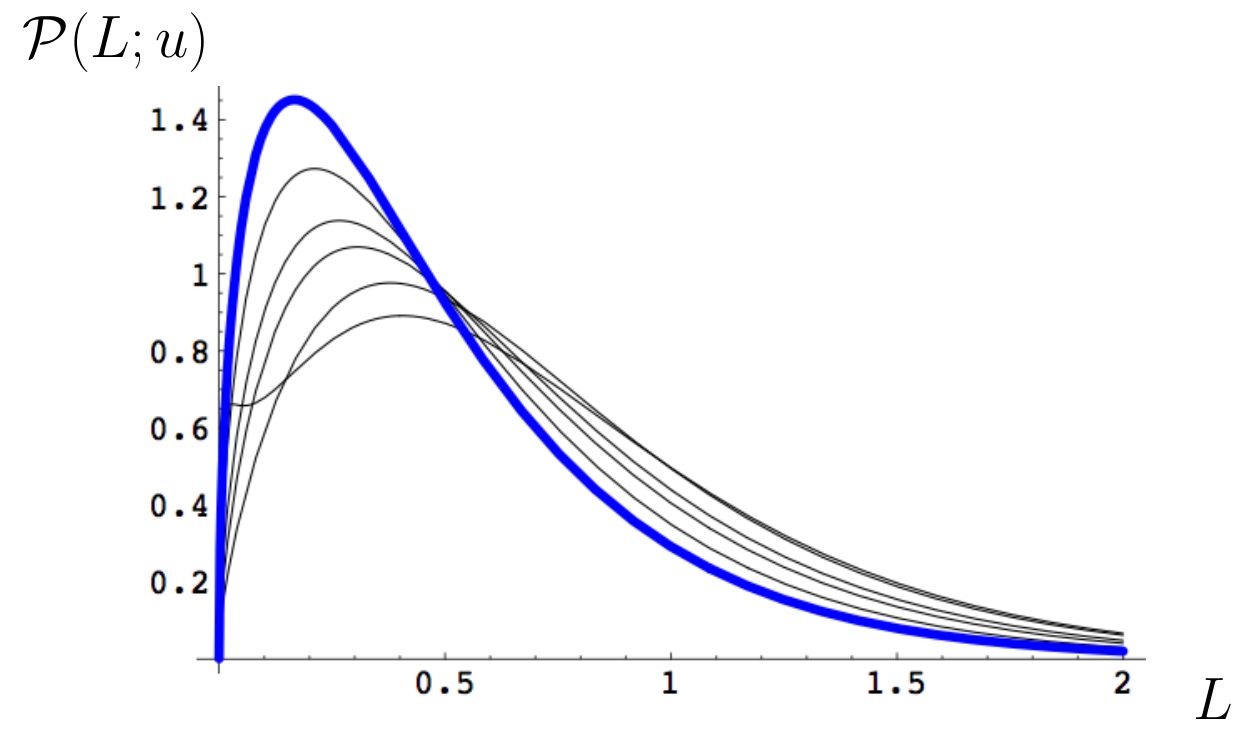}
\end{center}
\caption{The probability density $\mathcal{P}(L;u)$ for $c=1/3$ and for $u=0$ (blue thick line) and increasing values $u=1/8$, $1/4$, $1/3$, $1/2$ and $2/3$.
For this latter value, a peak starts to emerge around $L=0$.}
\label{fig:Proba1}
\end{figure} 

\subsubsection{Probability density for $L(d)$ when $d$ is a finite fraction of $k$.}
Our second result concerns the statistics of the hull perimeter at a distance $d$ corresponding to a \emph{finite
fraction} of the total distance $k$ between $v_0$ and $v_1$, in the limit of large $k$. In other words, we consider the situation where
\begin{equation*}
d= k\, u\ , \qquad 0 < u <1\ ,
\end{equation*}
for some fixed $u$ and for large $k$.
We find that:
\begin{equation}
\begin{split}
& \hskip -1.cm \lim_{k\to \infty} E_k (e^{-\tau L(k\, u)})= F\left(\sigma(\tau;u);u\right)
\quad \hbox{with}\  \sigma(\tau;u)\equiv 
\frac{1-2u+c\,  \tau  (1-u)^2}{u^2}\ ,\\
& \hskip -1.cmF(\sigma;u)\equiv \left(\frac{8 \sigma ^4+47 \sigma ^3+90 \sigma ^2+120 \sigma +48}{4 (\sigma +1)^{5/2}}-12 \right) \frac{u^3}{\sigma^3}\\ 
&\hskip -1.2cm\qquad \qquad  + \left(\frac{4 \sigma ^4+19 \sigma ^3+30 \sigma ^2+68 \sigma
   +32}{4 (\sigma +1)^{5/2}}+3 \sigma -8\right)\frac{u(1-2u)}{\sigma ^3}\\
& \hskip -1.2cm\qquad \qquad  +\left(\frac{4 \sigma ^5+16 \sigma ^4-7 \sigma ^2-40 \sigma -24}{4 (\sigma +1)^{5/2}}+3 \sigma ^2-5 \sigma +6\right) \frac{(1-2 u) (1-u)^2 \left(1-2u+2 u^2\right)}{\sigma ^4 u^3}\ ,\\
\label{eq:Ktau}
\end{split}
\end{equation}
with $c$ as in \eqref{eq:valc} above.

Note that, for $\tau\geq 0$ and $0<u<1$, $\sigma(\tau;u)+1=(1+c\, \tau)(1-u)^2/u^2>0$ so that the denominator $(\sigma+1)^{5/2}$ in
\eqref{eq:Ktau} does not vanish. For $0< u<1/2$, we have the stronger property $\sigma(\tau;u)>0$ and $F\left(\sigma(\tau;u);u\right)$ 
therefore clearly has no singularity for $\tau\geq 0$.
For $1/2<u< 1$ however, $\sigma(\tau;u)$ vanishes at the non-negative value 
$\tau=\frac{2u-1}{c\, (1-u)^2}$
and $F\left(\sigma(\tau;u);u\right)$ may seem at a first glance to develop some singularity there. Such behavior is not allowed for the Laplace 
transform of some probability density so a closer look at the formula is required. 
Expanding $F(\sigma;u)$ around $\sigma=0$ shows that $F(\sigma;u)$ is in fact well-behaved around $\sigma=0$ 
with $F(\sigma;u)=(512 u^6-3012 u^5+7518 u^4-10020 u^3+7515 u^2-3006 u+501)/(64 u^3)+O(\sigma)$. The seeming singularity 
at $\tau=\frac{2u-1}{c\, (1-u)^2}$ is therefore only an illusion.

Via an inverse Laplace transform, whose details are discussed in Appendix A, eq.~\eqref{eq:Ktau} is equivalent to\footnote{Here
$\text{erf}(a)$ denotes the usual error function $(2/\sqrt{\pi})\int_0^a e^{-z^2}\, dz$.}
\begin{equation}
\begin{split}
&  \lim_{k\to \infty} P_k (L\leq L(k\, u)< L+dL)=\mathcal{P}(L;u)\, dL \\
& \hbox{with}\ \mathcal{P}(L;u)= \frac{2}{\sqrt{\pi}}\frac{\sqrt{L}}{c^{3/2}}e^{-\frac{L}{c}}\times  \frac{\left(e^\frac{L}{c\, b}   \sqrt{\frac{\pi L}{c\, b}  } \left(1-\text{erf}\left(\sqrt{\frac{L}{c\, b}  }\right)\right)-1\right) p\left(\frac{L}{c\, b}\right)+r\left(\frac{L}{c\, b}\right)}{4\, 
(\sqrt{b}+b)^3 }\\
& \\
   & {\rm where}\ b\equiv\ b(u)=\frac{(1-u)^2}{u^2}\ , \\
  & p(\ell )\equiv 2 b \left(b^2-1\right) \ell ^2-\left(5 b^3+3 b+4\right) \ell +6 \left(b^3-1\right)\ , \\
  &  r(\ell )\equiv  b \left(15 b^2-1\right) \ell+2 \left(5 b^3-1\right)\ .\\
\end{split}
\label{eq:Pk}
\end{equation}
The function $\mathcal{P}(L;u)$ is plotted in figure \ref{fig:Proba1} for various values of $u$ and $c=1/3$. The limit $u\to 0$ (i.e $b\to \infty$)
describes situations where the distance $d$ at which the hull perimeter is measured does not scale with $k$, 
so we expect to recover the result of eq.~\eqref{eq:Pinf} in this limit. This is indeed the case as: 
\begin{equation*}
\mathcal{P}(L;0)= \frac{2}{\sqrt{\pi}}\frac{\sqrt{L}}{c^{3/2}}e^{-\frac{L}{c}}\ .
\end{equation*}
The probability density $\mathcal{P}(L;0)$ is displayed in blue in figure \ref{fig:Proba1} (thick line). Note that $\mathcal{P}(L;u)/\mathcal{P}(L;0)$ is a function
of the variable
\begin{equation*}
R\equiv \frac{L}{b(u)}=L\, \frac{u^2}{(1-u)^2}= L\, \frac{d^2}{(k-d)^2}\ .
\end{equation*}
Recall that $L(d)$ is the ratio of the the actual hull perimeter $\mathcal{L}(d)$ by its ``natural scale" $d^2$. The new variable 
$R$ therefore corresponds to probing the value of the random variable
\begin{equation*}
R(d)\equiv \frac{\mathcal{L}(d)}{(k-d)^2}\ ,
\end{equation*}
i.e.\ measuring the hull perimeter at a scale no longer fixed by $d^2$ by rather by $(k-d)^2$.
\begin{figure}
\begin{center}
\includegraphics[width=9cm]{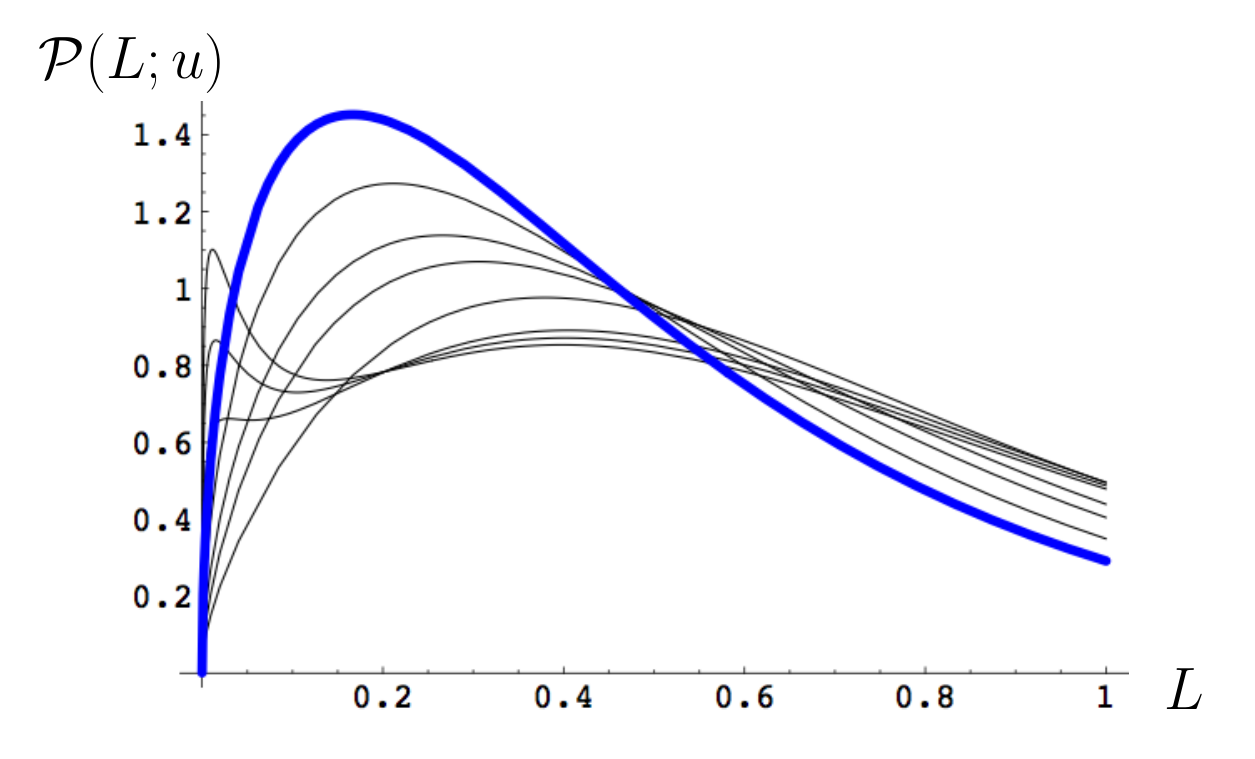}
\end{center}
\caption{The same plot as figure ~\ref{fig:Proba1} limited to the range $0\leq L\leq 1$, with two new value $u\simeq 0.69$ and $u\simeq 0.71$.
This plot is to emphasize the emergence of a peak for small $L$ when $u\geq 1/2$.}
\label{fig:Proba2}
\end{figure}
When $u$ becomes larger than $1/2$, a peak for small $L$ starts to emerge in $\mathcal{P}(L,u)$, as displayed in figure \ref{fig:Proba2}. 
This peak increases when $u$ tends to $1$ ($b\to 0$). This limit is best captured by switching to the variable $R$, namely considering 
\begin{equation*}
\lim_{k\to \infty} P_k (R\leq R(k\, u)< R+dR)= \tilde{\mathcal{P}}(R;u)\, dR\quad \hbox{with}\quad \tilde{\mathcal{P}}(R;u)= b\, \mathcal{P}(b\, R;u) \ .
\end{equation*}
\begin{figure}
\begin{center}
\includegraphics[width=8.5cm]{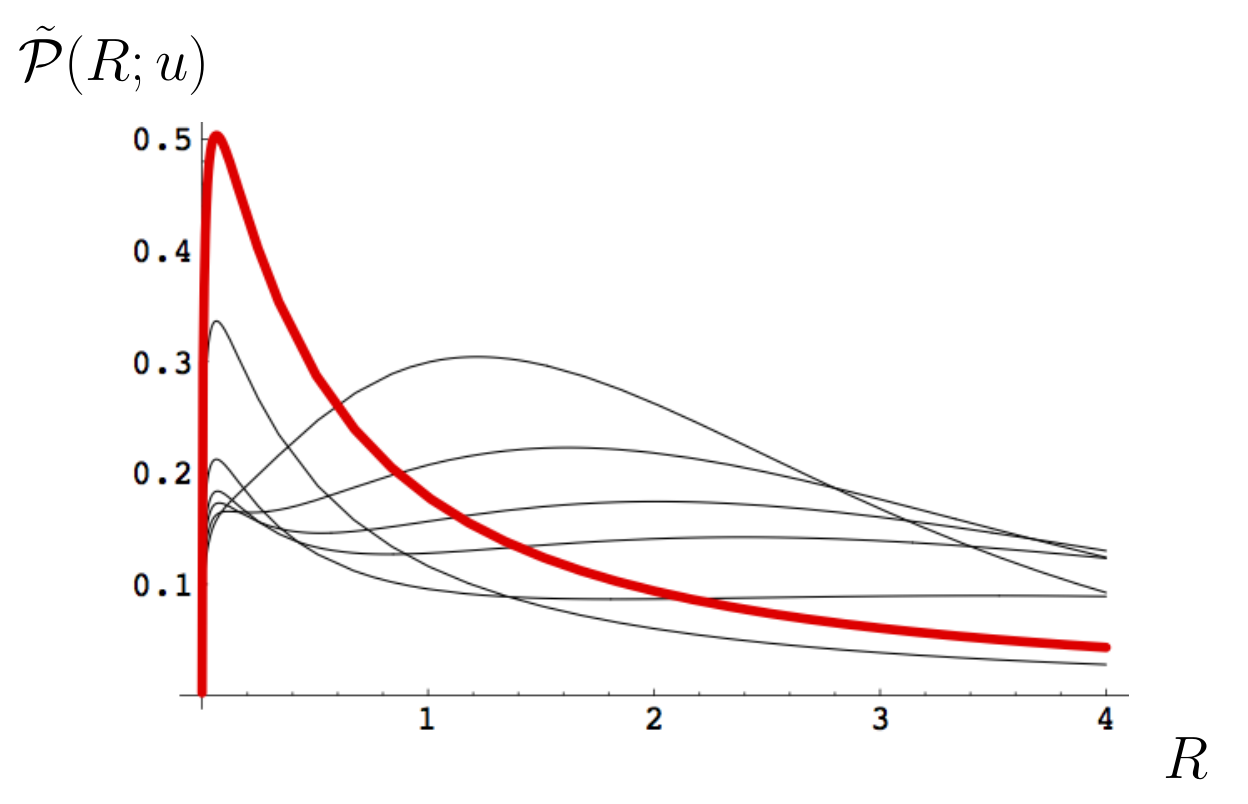}
\end{center}
\caption{The probability density $\tilde{\mathcal{P}}(R;u)$ for $c=1/3$ and for $u\simeq 0.64$, $0.67$, $0.69$, $0.71$, $0.75$, $0.875$
and $u=1$ (red thick line). The height of the peak for small $L$ increases with increasing $u$.}
\label{fig:Proba3}
\end{figure}
In other words, when $d$ approaches $k$, the natural scale for $\mathcal{L}(d)$ is no longer $d^2$ but rather $(k-d)^2$.
The probability density $\tilde{\mathcal{P}}(R,u)$ is displayed in figure \ref{fig:Proba3} for $c=1/3$ and various values of $u$. 
When $u\to 1$ ($b\to 0$), $\tilde{\mathcal{P}}(R;u)$ 
converges to a well-defined distribution
\begin{equation*}
\tilde{\mathcal{P}}(R;1)=2 \sqrt{\frac{R}{\pi c^5}} (R+c)-\frac{R}{c^3}\, (2 R+3c)\, 
  e^\frac{R}{c} \left(1-\text{erf}\left(\sqrt{\frac{R}{c}}\right)\right)\ .
\end{equation*}
\begin{figure}
\begin{center}
\includegraphics[width=8.5cm]{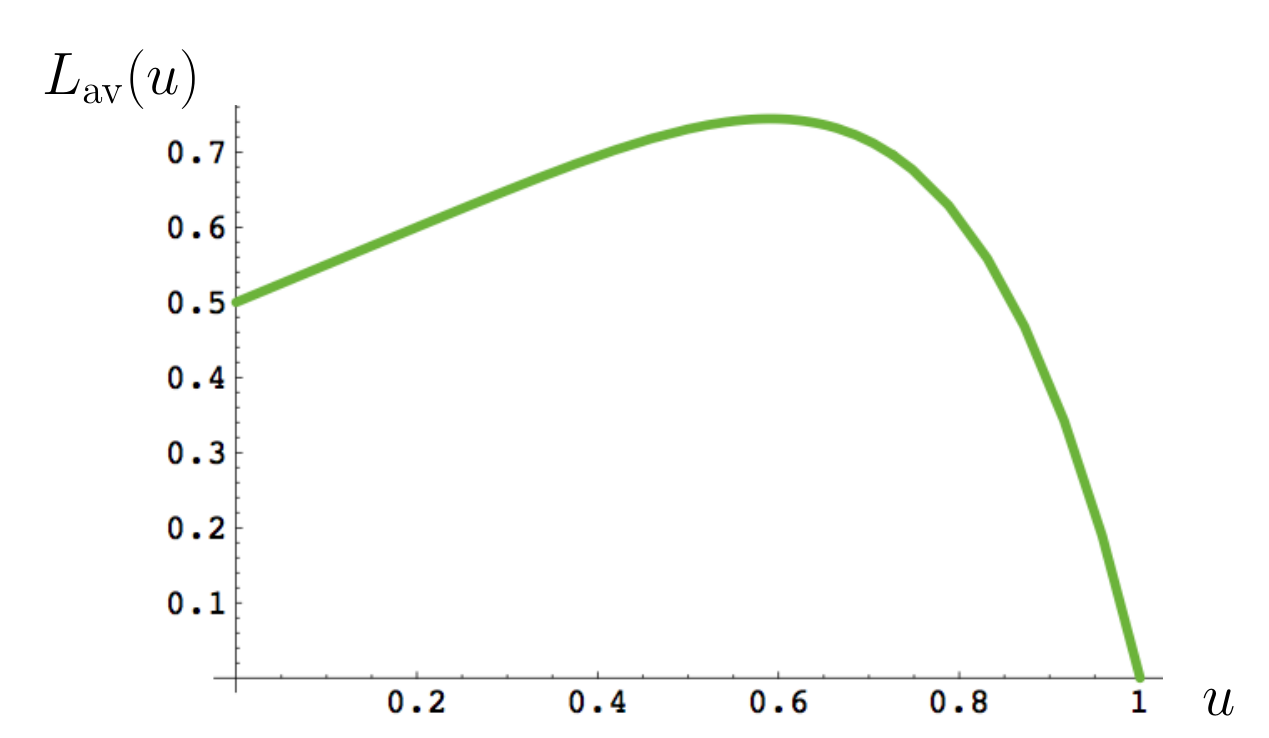}
\end{center}
\caption{The average value $L_{\rm av}(u)\equiv  \lim_{k\to \infty} E_k (L(k\, u))$ as a function of $u$ (here $c=1/3$).}
\label{fig:Av1}
\end{figure}
\begin{figure}
\begin{center}
\includegraphics[width=10cm]{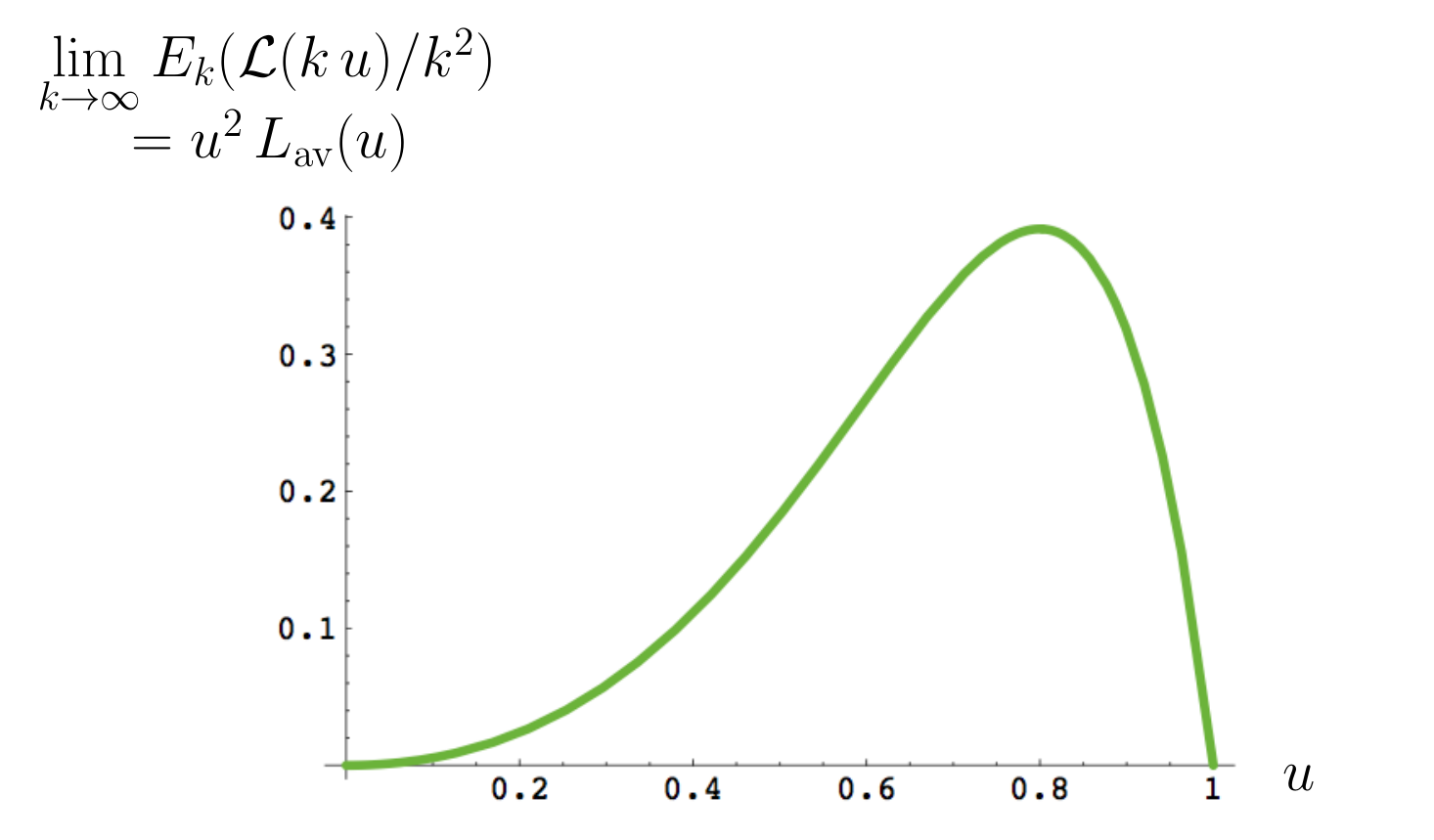}
\end{center}
\caption{The average profile of the hull perimeter, i.e.\ the quantity $\lim_{k\to \infty} E_k (\mathcal{L}(k\, u)/k^2)= u^2\, L_{\rm av}(u)$ as a function of $u$ (here $c=1/3$).}
\label{fig:Av2}
\end{figure}
This distribution is displayed in red in figure \ref{fig:Proba3} (thick line). Note that this distribution has all its moments \emph{infinite}. 

This result may appear strange at a first glance but the divergence of, say the first moment is in fact consistent with a direct computation of the average values 
of $L(d)$ and $R(d)$ for arbitrary $u$:
expanding eq.~\eqref{eq:Ktau} at first order in $\tau$, we have indeed
\begin{equation*}
L_{\rm av}(u)\equiv \lim_{k\to \infty} E_k (L(k\, u)) = \frac{3c}{2}(1+u-3 u^6
+u^7) \\
\end{equation*}
The average value $L_{\rm av}(u)$ is displayed in figure \ref{fig:Av1} for $u<0<1$ and $c=1/3$. When $u\to 0$, it tends to a finite value
$3c/2$, meaning that the average hull perimeter at distance $d$ from the origin scales like $3c/2 \times d^2$ 
when $k$ is infinitely large, as expected. When $d$ corresponds to a finite fraction $u<1$ of $k$,  the average hull perimeter remains of order $d^2$
with a finite prefactor $L_{\rm av}(u)$ depending on $u$. When $d\to k$ however, i.e. $u\to 1$, then $L_{\rm av}(u)\sim 15c(1-u)$ hence 
the average hull perimeter vanishes like $15 c\, k(k-d)$. 
This vanishing is not surprising since the hull perimeter is also the length of the boundary of the domain $\mathcal{C}_d$ in which the vertex 
$v_1$ is ``trapped". When $d$ approaches $k$, this domain becomes smaller and smaller and so does its boundary. This 
vanishing is only linear in $(k-d)$ and the average value of $R(d)$ behaves accordingly as $15 c\, k^3/(k-d)$, hence diverges when $d\to k$.

We may finally consider the average ``profile" of the hull perimeter, i.e.\ the average value $E_k (\mathcal{L}(k\, u)/k^2)$
of the hull perimeter at distance $d=k\, u$ normalized for all $u$ by the same \emph{global scale} $k^2$ (instead of the local natural scale $(k\, u)^2$).
In the limit $k\to \infty$, it is simply equal to $u^2\, L_{\rm av}(u)$ and has the form displayed in figure \eqref{fig:Av2}.

\begin{figure}
\begin{center}
\includegraphics[width=8cm]{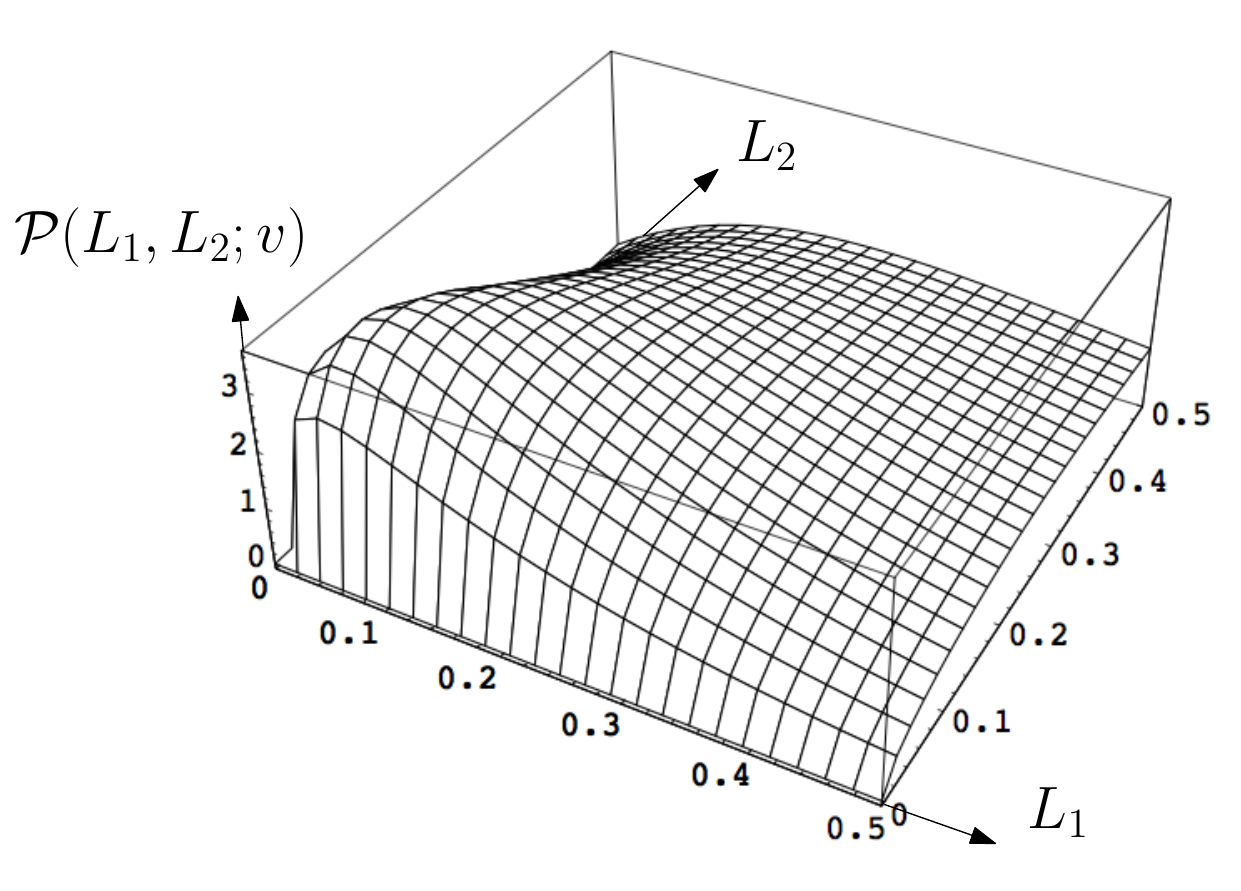}
\end{center}
\caption{The joint probability density $\mathcal{P}(L_1,L_2;v)$ for $v=2$ (here $c=1/3$).}
\label{fig:ProbL1L2}
\end{figure} 

\subsubsection{Joint probability density for $L(d_1)$ and $L(d_2)$ when $k\to \infty$}
Our third main result deals with the joint law for hull perimeters at distances $d_1$ and $d_2$ from $v_0$, with again $k\to \infty$.
Assuming $d_2>d_1$, we set:
 \begin{equation*}
 d_1=d\ ,  \quad d_2=v\, d\ , \qquad\hbox{with}\  v>1\ .
 \end{equation*}
 We find that:
\begin{equation}
\begin{split}
& \hskip -1.cm \lim_{d\to \infty } E_\infty (e^{-\tau_1 L(d)-\tau_2 L(v d)})\\
&  \hskip -1.2cm \qquad = \frac{v^3}{\left(v^2 (1\!+\!c\, \tau_1) (1\!+\!c\,  \tau_2)\!-\!2\, v\, c\, \tau_2  \left(1\!+\!c\, \tau_1\! -\!\sqrt{1\!+\!c\, \tau_1
   }\right)\!+\!c\, \tau_2  \left(2\!+\!c\, \tau_1\! -\!2 \sqrt{1\!+\!c
   \, \tau_1}\right)\right)^{3/2}} \\
   \label{eq:zuniv}
\end{split}
\end{equation}
which equivalently yields a joint probability density (see Appendix A for details on the appropriate double inverse Laplace transform):
\begin{equation}
\begin{split}
& \hskip -.7cm \lim_{d\to \infty} P_\infty \Big(L_1\leq L(d)< L_1+dL_1\ \hbox{and}\  L_2\leq L(v\, d)< L_2+dL_2\Big)=\mathcal{P}(L_1,L_2;v)\, dL_1\, dL_2 \ ,\\
& \mathcal{P}(L_1,L_2;v)= \frac{2}{\sqrt{\pi}}  \frac{\sqrt{L_1}}{c^{3/2}} e^{-\frac{L_1}{c}} \times \frac{\sqrt{2}}{c}\frac{v^2 e^{-\frac{L_2\, v^2}{c\, (v-1)^2}}}{(v-1)^2} 
\sum _{n=0}^\infty \frac{(-1)^n  \left(\frac{L_1}{c (v-1)^2}\right)^{\frac{n}{2}}
\pi_n\left(\sqrt{\frac{2 L_2\, v^2}{c\, (v-1)^2}}\right)}{(n+1)!\, \Gamma \left(\frac{n+1}{2}\right)}\ . \\ 
\end{split}
\label{eq:PL1L2}
\end{equation}
Here $\pi_n(t)$ is the polynomial (of the Hermite type) defined by\footnote{The first polynomials are
$\pi_0(t)=t$, $\pi_1(t)=t(3-t^2)$, $\pi_2(t)=t(12-9t^2+t^4)$, $\pi_3(t)=t(60-75t^2+18 t^4-t^6)$.}
\begin{equation}
\pi_n(t)\equiv - e^{\frac{t^2}{2}} \frac{d^{n+1}}{dt^{n+1}}\left(t^n\, e^{-\frac{t^2}{2}}\right)\ .
\label{eq:pin}
\end{equation}
The joint probability density $\mathcal{P}(L_1,L_2;v)$ is plotted for $v=2$ in figure \ref{fig:ProbL1L2}.

\begin{figure}
\begin{center}
\includegraphics[width=8cm]{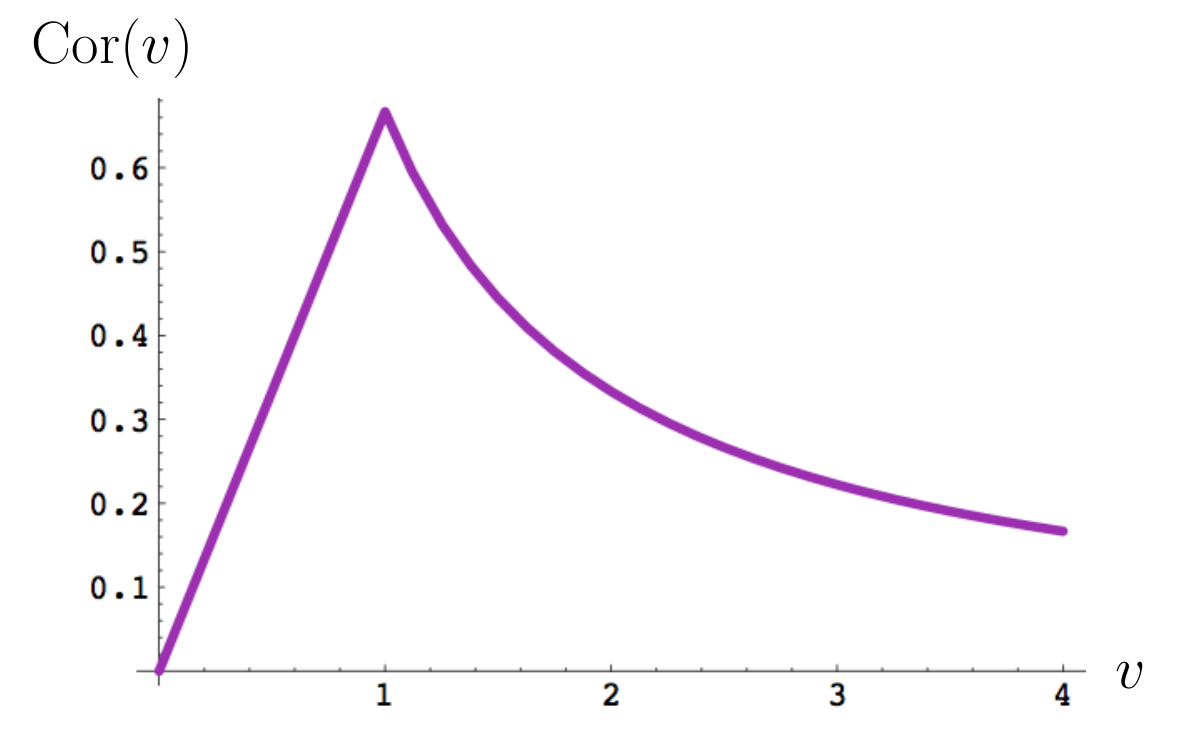}
\end{center}
\caption{The correlation ${\rm Cor}(v)$ as a function of $v$. This quantity is independent of $c$.}
\label{fig:Cor}
\end{figure} 
Expanding eq.~\eqref{eq:zuniv} at first order in $\tau_1$ and $\tau_2$, we immediately get
\begin{equation*}
\lim_{d \to \infty} E_\infty ( L(d)\, L(v \, d))= \frac{3}{4} c^2 \left(3+\frac{2}{v}\right)\ \hbox{for}\ v>1
\end{equation*}
hence, by dividing by the common average value $3c/2$ for $L(d)$ and $L(v\, d)$ (and using an obvious symmetry to extend the result to $v<1$), a correlation
 \begin{equation*}
\hskip -1.cm {\rm Cor}(v)\equiv \frac{ \lim_{d \to \infty} E_\infty ( L(d)\, L(v \, d))}{\lim_{d \to \infty} E_\infty ( L(d))\times \lim_{d \to \infty} E_\infty (L(v\, d))}-1
=\frac{2}{3\max(v,1/v)}\qquad  \hbox{for all}\ v>0\ .
\end{equation*}
Note that this correlation is independent of $c$. The value $2/3$ at $v=1$ is obtained directly from eq.~\eqref{eq:CLG}.
The correlation ${\rm Cor}(v)$ is plotted in figure \ref{fig:Cor} for illustration. 

Another measure of the correlation is the average value $L_{\rm av}(v|L_1)$ 
of $L(d\ v)$, knowing that $L(d)$ lies in the range $[L_1,L_1+dL_1]$. By integrating $L_2\, \mathcal{P}(L_1,L_2;v)$ over
$L_2$, it is easily found to be\footnote{
It is easily verified that, when computing the $p$-th moment of $L_2$ with the distribution $\mathcal{P}(L_1,L_2;v)$, only the first $2p+1$ polynomials
$\pi_n$ (i.e.\ $n=0,\cdots ,2p$) in \eqref{eq:PL1L2} give a non-zero contribution.}:
\begin{equation*}
L_{\rm av}(v|L_1)= \frac{3 c\, (v-1)^2+3 \sqrt{\pi\, c\,  L_1}\, (v-1)+2\, L_1}{2 v^2}\ , \qquad v>1\ ,
\end{equation*}
which varies from $L_{\rm av}(1|L_1)=L_1$ to $L_{\rm av}(\infty|L_1)=3c/2$, as expected.

\section{The slice recursion and the hull perimeter}
\label{sec:recursion}
We now come to the derivation of our various results. It relies on the existence of a recursion
relation for slice generating functions, as described in \cite{G15b} for quadrangulatiions and \cite{G15a} for triangulations.
As we shall see, the origin of this recursion is indeed intimately linked to the notion of hull boundary and this will eventually 
allows us to have a direct control on the hull perimeter.

\subsection{The case of quadrangulations}
\label{sec:quad}
\subsubsection{The slice recursion}
We consider here the $k$-slices defined in Section~\ref{sec:definitions} in correspondence with $k$-pointed-rooted quadrangulations,
and more generally the larger set of $\ell$-slices with $1\leq \ell \leq k$. Recall that, in an $\ell$-slice, $\ell$ is the length of
its left boundary and is also the distance $d(v_0,v_1)$ between its apex and the first extremity of its base.
Let us denote by $R^{(q)}_k\equiv R^{(q)}_k(g)$, $k\geq 1$, the generating function for this larger family of slices,
where we assign a weight $g$ to each tetravalent \emph{inner} face (i.e.\ each face other than the outer face). 
The quantity $R^{(q)}_k$ is also the generating function, with a weight $g$ per face, of pointed-rooted quadrangulations\footnote{By 
pointed-rooted quadrangulations, we mean in general $\ell$-pointed-rooted quadrangulations with arbitrary $\ell\geq 1$.} whose graph distance $d(v_0,v_1)$ between 
the origin $v_0$ and the first extremity $v_1$ of the root edge satisfies $1\leq d(v_0,v_1)\leq k$.
The explicit expression for $R^{(q)}_k$ can be found in \cite{G15b} and reads:
\begin{equation*}
R^{(q)}_k=R^{(q)} \frac{(1-x^{k})(1-x^{k+3})}{(1-x^{k+1})(1-x^{k+2})}\ , \qquad R^{(q)}=\frac{1-\sqrt{1-12g}}{6g}\ ,
\end{equation*}
where $x$ is a parametrization of $g$ through
\begin{equation*}
g=\frac{x(1+x+x^2)}{(1+4x+x^2)^2}\ .
\end{equation*}
Note that $x$ and $1/x$ lead to the same value of $g$ so we shall impose the extra condition $|x|\leq 1$ to univocally fix $x$.
The generating functions are well-defined for real $g$ in the range $0\leq g\leq 1/12$: this then corresponds to a real $x$ in the range 
$0\leq x\leq 1$.
To be precise, the recursion relation found in \cite{G15b} concerns the quantity
\begin{equation*}
T^{(q)}_k\equiv R^{(q)}_k-R^{(q)}_1\ , \qquad k\geq 1\ ,
\end{equation*}
which enumerates $\ell$-slices whose left boundary length $\ell$ is between $2$ and $k$ (note that $T^{(q)}_1=0$). 
From the explicit expression of $R^{(q)}_k$, we immediately deduce
\begin{equation*}
T^{(q)}_k=T^{(q)} \frac{(1-x^{k-1})(1-x^{k+4})}{(1-x^{k+1})(1-x^{k+2})}\ , \qquad T^{(q)}= \frac{x(1+4x+x^2)}{(1+x+x^2)^2}\ .
\end{equation*}
It was shown in \cite{G15b} that this generating function satisfies a recursion relation of the form: 
\begin{equation}
\begin{split}
& T^{(q)}_k=\mathcal{K}^{(q)}(T^{(q)}_{k-1})\ ,\\ 
& \mathcal{K}^{(q)}(T)=\frac{(R^{(q)}_1)^2 (T+R^{(q)}_1)\, \Phi^{(q)}(T)}{1-R^{(q)}_1(T+R^{(q)}_1)\, \Phi^{(q)}(T)}\ , \  \Phi^{(q)}(T)\equiv \Phi^{(q)}(T,g)=\sum_{i\geq 2} h^{(q)}_{2i}(g)\, T^{i-2}\ ,\\
\end{split}
\label{eq:recrel}
\end{equation}
for $k\geq 2$ with $T^{(q)}_1=0$.
Here $h^{(q)}_{2i}(g)$ are appropriate generating functions whose definition can be found in \cite{G15b} and whose explicit 
expression will not be needed in our calculation. The precise form of the ``kernel" $\mathcal{K}^{(q)}(T)$ is also not important 
for our calculation and we displayed it only to help the reader make the connection with \cite{G15b}. What matters for us
is only the following simple property: the kernel $\mathcal{K}^{(q)}(T)$ is independent of $k$ whereas
$T^{(q)}_k$ depends on $k$ only through the variable $x^k$. 
We immediately deduce that, if we make the transformation $x^k\to \lambda x^k$ in the expressions
for $T^{(q)}_k$, the obtained quantity still satisfies the same recursion relation. In other words, if we set
\begin{equation}
T^{(q)}_k(\lambda)= T^{(q)}\frac{(1-\lambda\, x^{k-1})(1-\lambda\, x^{k+4})}{(1-\lambda\, x^{k+1})(1-\lambda\,  x^{k+2})}\ ,
\label{eq:Tlambdak}
\end{equation} 
then $T^{(q)}_k(\lambda)$ still satisfies\footnote{A detailed calculation shows that this actually holds only for $\lambda$ close enough to $1$ so that
we are guaranteed that  $\frac{(1-\lambda^2\, x^{2 k+1})}{\left(1-\lambda\, x^{k}  \right) \left(1-\lambda\, x^{k+1}\right)}$ remains non-negative.
We shall be in this regime in the following.}
\begin{equation}
T^{(q)}_k(\lambda)=\mathcal{K}^{(q)}(T^{(q)}_{k-1}(\lambda))\ .
\label{eq:reclambda}
\end{equation}
This, together with the explicit form \eqref{eq:Tlambdak} of $T^{(q)}_k(\lambda)$ is the only ingredient that we shall rely on in the following for our explicit calculations.

\subsubsection{Connection with the hull perimeter}
\begin{figure}
\begin{center}
\includegraphics[width=9cm]{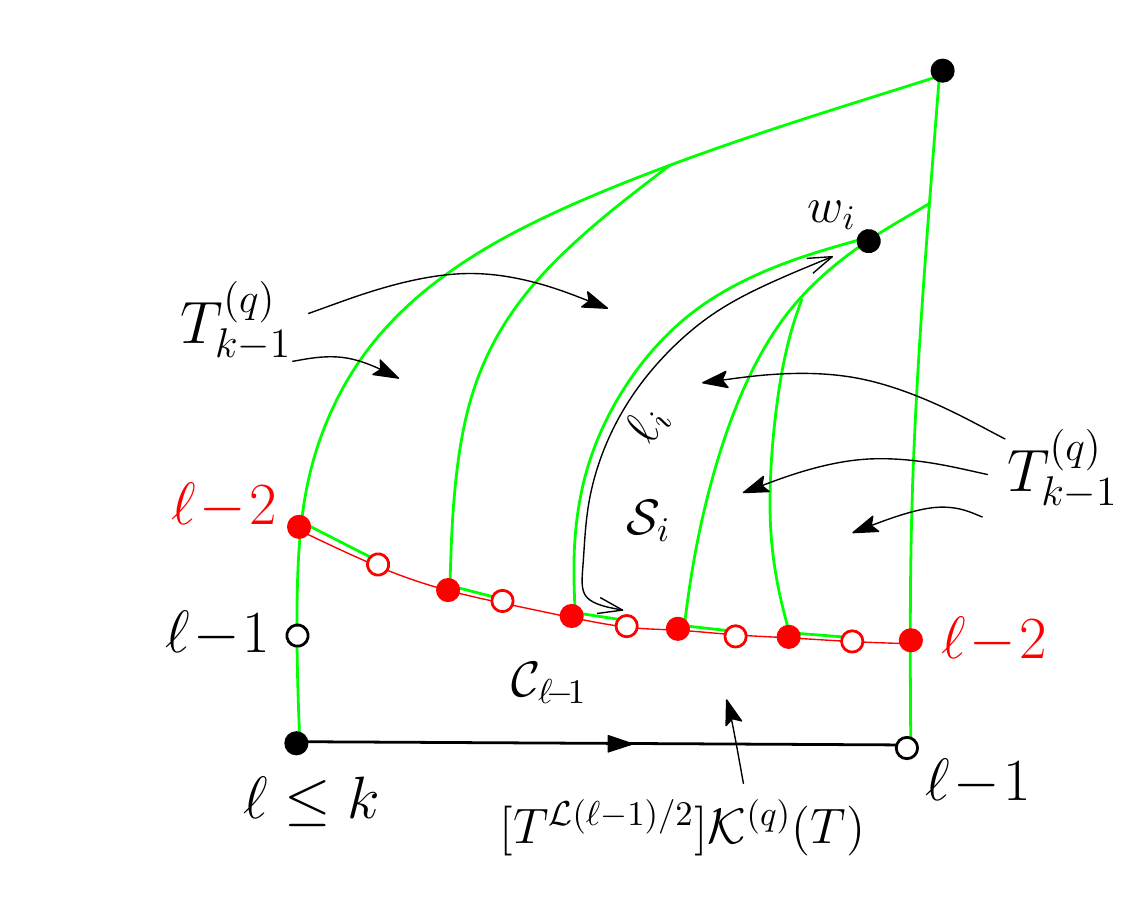}
\end{center}
\caption{The decomposition of a slice enumerated by $T^{(q)}_k$ leading to the recursion relation \eqref{eq:recrel}. The hull boundary
at distance $\ell-1$, represented in red and made of an alternating sequence of vertices at distance $\ell-2$ (filled red circles)
and at distance $\ell-1$ (open red circles) from the apex, serves as a dividing line of even length $\mathcal{L}(\ell-1)$ (here equal to $10$)
separating a domain $\mathcal{C}_{\ell-1}$ from the
hull. By drawing the leftmost shortest paths to the apex from each of the vertices at distance $\ell-1$ along the hull boundary (open red circles), the hull itself is decomposed into a number $\mathcal{L}(\ell-1)/2$ of sub-slices, each enumerated by $T^{(q)}_{k-1}$
(since each sub-slice $\mathcal{S}_i$ has a left boundary length $\ell_i$ satisfying $2\leq \ell_i\leq k-1$). The generating function
of the domain $\mathcal{C}_{\ell-1}$ is $[T^{\mathcal{L}(\ell-1)/2}]\mathcal{K}^{(q)}(T)$.}
\label{fig:constrhull}
\end{figure} 
Let us now briefly recall the origin of the recursion relation and show how it may allow us to control the hull perimeter. 
Starting with an $\ell$-slice with left boundary length $\ell$ between $2$ and $k$ (as enumerated by $T^{(q)}_k$),
the recursion is obtained by cutting the $\ell$-slice along some particular line, called the ``dividing line" in \cite{G15b}. This dividing line is 
precisely the hull boundary \emph{at distance $\ell-1=d(v_0,v_1)-1$} from $v_0$ in the $\ell$-slice, as we defined it in 
Section~\ref{sec:definitions}\footnote{In \cite{G15b}, a first edge of the right boundary of the slice, linking its vertices $v_0^{(\ell)}$ and $v_0^{(\ell-1)}$ at respective distances 
$\ell-1$ and $\ell-2$ from $v_0$ was added for convenience to the dividing line. This edge is not present in our definition where we let the dividing line
start at the vertex $v_0^{(\ell-1)}$ at distance $\ell-2$.} (see figure \ref{fig:constrhull}). The generating
function $T^{(q)}_k$ is obtained by multiplying the generating function of the domain $\mathcal{C}_{\ell-1}$ (corresponding in the slice to the domain 
on the same side of the hull boundary as $v_1$) by the generating function of the domain $\mathcal{H}_{\ell-1}$ (corresponding in the slice to the domain 
on the same side of the boundary as $v_0$). For a \emph{fixed value} $\mathcal{L}$ of the hull perimeter 
$\mathcal{L}(\ell-1)$, the first generating function is easily seen to be independent of $\ell$ (since there is no restriction on mutual distances
within this domain apart from the fact that $v_1$ is at distance $1$ from the boundary). 
This generating function is nothing but (see \cite{G15b} for details):
\begin{equation*}
[T^{\mathcal{L}/2}]\mathcal{K}^{(q)}(T)
\end{equation*}
(recall that $\mathcal{L}$ is even).
As for the domain $\mathcal{H}_{\ell-1}$, it is formed of \emph{exactly $\mathcal{L}/2$ slices} with respective left boundary lengths $\ell_1,\ell_2,\cdots,\ell_{\mathcal{L}/2}$, each satisfying $2\leq \ell_i\leq \ell-1$, hence $2\leq \ell_i\leq k-1$ when we sum over all possible values of $\ell$ between 
$2$ and $k$. These slices are obtained by decomposing the domain $\mathcal{H}_{\ell-1}$ upon cutting along the leftmost shortest paths to $v_0$
from the ${\mathcal{L}/2}-1$ vertices of its boundary which are at distance $\ell-1$ from $v_0$, but the last one\footnote{The leftmost shortest path
from this last vertex to $v_0$ follows the left boundary of the original $\ell$-slice and need not being cut.}.
Note that when $\ell=2$, the hull boundary at distance $\ell-1=1$ must be understood as reduced to the single
vertex $v_0$, having length $0$. Each of the $\mathcal{L}/2$ slices composing the domain $\mathcal{H}_{\ell-1}$ is thus enumerated by 
$T^{(q)}_{k-1}$ and the net contribution of this domain to $T^{(q)}_k$ is eventually
\begin{equation*}
\left(T^{(q)}_{k-1}\right)^{\mathcal{L}/2}\ .
\end{equation*}
Summing over all (even) values of $\mathcal{L}$ yields the desired recursion relation \eqref{eq:recrel}. 

If we now wish to keep a control on the 
value $\mathcal{L}$ of the hull perimeter at distance $d(v_0,v_1)-1$ in the original $\ell$-slices (characterized by $2\leq d(v_0,v_1)\leq k$)
by assigning, say a weight $\alpha^{\mathcal{L}}$ to these $\ell$-slices, we simply need to replace $T^{(q)}_{k-1}$ by 
$\alpha^2\, T^{(q)}_{k-1}$ at the step $k-1\to k$ of the recursion. 
In other words, the generating function of $\ell$-slices with $2\leq \ell\leq k$, with a weight $g$ per inner face and a weight 
$\alpha^{\mathcal{L}(\ell-1)}$ is simply given by
\begin{equation*}
\sum_{\mathcal{L}\geq 0\atop \mathcal{L}\ {\rm even}}[T^{\mathcal{L}/2}]\mathcal{K}^{(q)}(T) \times \alpha^{\mathcal{L}} \times \left(T^{(q)}_{k-1}\right)^{\mathcal{L}/2}=
 \mathcal{K}^{(q)}(\alpha^2\, T^{(q)}_{k-1})\ .
\end{equation*}
In the argument leading to this formula, $\ell$-slices with $\ell=2$ must be understood as having $\mathcal{L}(\ell-1)=0$ (since the hull boundary is reduced to $v_0$ in this case), in agreement with our general convention $\mathcal{L}(1)=0$ for slices associated with quadrangulations.
The $2$-slices are thus enumerated with a weight $\alpha^{\mathcal{L}(1)}=1$.
\begin{figure}
\begin{center}
\includegraphics[width=8cm]{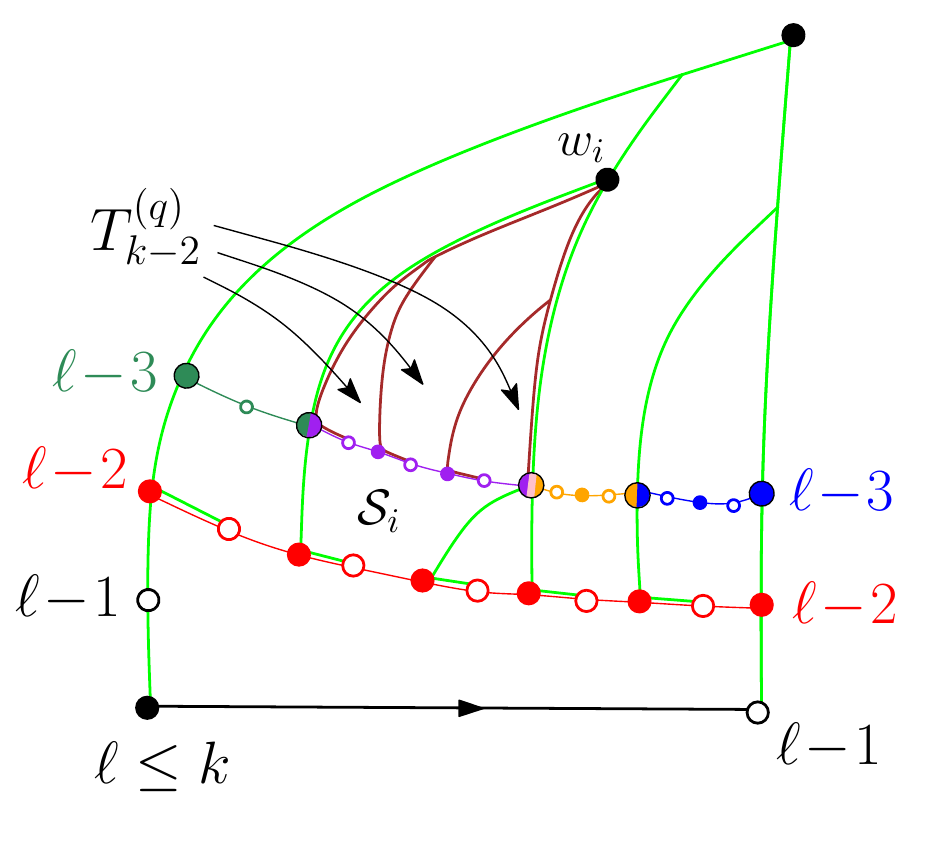}
\end{center}
\caption{The construction of the hull boundary at distance $\ell-2$ (i.e.\ formed by an alternating sequence of vertices at distance
$\ell-3$ and $\ell-2$) by concatenating the hull boundaries at distance $\ell_i-1$ (represented in different colors) 
of the various sub-slices $\mathcal{S}_i$ forming $\mathcal{H}_{\ell-1}$.}
\label{fig:constrhull2}
\end{figure} 

\subsubsection{Controlling the hull perimeter at some arbitrary $d$}
If we wish instead to control the hull perimeter at distance $\ell-2$ from the apex $v_0$ in $\ell$-slices with $3\leq \ell \leq k$,
 we may simply repeat our construction within each of the $\mathcal{L}(\ell-1)/2$ sub-slices $\mathcal{S}_i$ forming the domain $\mathcal{H}_{\ell-1}$
(recall that the left boundary length $\ell_i$ of the sub-slice $\mathcal{S}_i$ satisfies $2\leq \ell_i\leq k-1$). More precisely we start by constructing
the hull boundary at distance $\ell_i-1$ within each sub-slice $\mathcal{S}_i$. Here the distances within the sub-slice $\mathcal{S}_i$ are measured
\emph{from its apex} $w_i$ which serves as origin of the sub-slice. In particular, if $\ell_i=2$ for some $i$, its hull boundary is reduced to the vertex
$w_i$. The hull boundary at distance $\ell-2$ is then obtained by \emph{concatenating} all the hull boundaries
at distance $\ell_i-1$ of the successive sub-slices $\mathcal{S}_i$ (see figure \ref{fig:constrhull2})\footnote{Note that when $\ell_i=2$ for some $i$, the hull boundary at distance $\ell-2$ simply passes through the apex $w_i$ without entering the sub-slice $\mathcal{S}_i$, hence contributes $0$ to the hull perimeter.}. This property is a direct consequence of the fact that the notions of distances and leftmost shortest paths within the $\ell$-slice 
are strictly bound to the same notions within the sub-slices $\mathcal{S}_i$. In particular, even though 
the apex $w_i$ of $\mathcal{S}_i$ is in general distinct from $v_0$, the 
distance $d(v_0,v)$ from $v_0$ to any vertex $v$ inside $\mathcal{S}_i$ is equal to $d(v_0,w_i)$ plus the distance \emph{within the
slice $\mathcal{S}_i$} from $w_i$ to $v$: this is a direct consequence of the fact that the sub-slice boundaries are shortest paths from
their base extremities to $v_0$ in the original $\ell$-slice. 

To summarize, the boundary perimeter $\mathcal{L}(\ell-2)$ is the sum of the hull perimeters at distance $\ell_i-1$ of the
$\mathcal{L}(\ell-1)/2$ sub-slices $\mathcal{S}_i$ forming the domain $\mathcal{H}_{\ell-1}$. As before, each of these sub-slice hull perimeters
corresponds to twice the number of sub-sub-slices forming the hull at distance $\ell_i-1$ of the sub-slice $\mathcal{S}_i$ at hand,
each sub-sub-slice being now enumerated by $T^{(q)}_{k-2}$ (see figure \ref{fig:constrhull2}). Considering two consecutive steps of our recursion relation, 
we deduce that the quantity
\begin{equation*}
 \mathcal{K}^{(q)}\left( \mathcal{K}^{(q)}\left(\alpha^2\, T^{(q)}_{k-2}\right)\right)\ -  \mathcal{K}^{(q)}(0).
\end{equation*}
is the generating function of $\ell$-slices with $3\leq \ell\leq k$ with a weight $g$ per inner face and a weight $\alpha^{\mathcal{L}(\ell-2)}$.
As before, we have taken the convention that $\ell$-slices with $\ell=3$ have perimeter $\mathcal{L}(\ell-2)=0$ (since their hull boundary is reduced
to $v_0$). The subtracted term $ \mathcal{K}^{(q)}(0)=T^{(q)}_1$ suppresses the $\ell$-slices with $\ell=2$ which would otherwise
be present from the first term. We may instead consider the un-subtracted quantity 
\begin{equation*}
 \mathcal{K}^{(q)}\left( \mathcal{K}^{(q)}\left(\alpha^2\, T^{(q)}_{k-2}\right)\right)
\end{equation*}
which is the generating function of $\ell$-slices with $2\leq \ell\leq k$ with a weight $g$ per inner face and a weight $\alpha^{\mathcal{L}(\ell-2)}$
\emph{if} $\ell > 2$ (and no $\alpha$-dependent weight otherwise) .
Repeating the argument, we may, for $1\leq m<k$, identify 
\begin{equation}
\underbrace{\mathcal{K}^{(q)}\big( \mathcal{K}^{(q)}\big(\cdots \big( \mathcal{K}^{(q)}\big(}_{m\ \hbox{\scriptsize times}}\alpha^2\, T^{(q)}_{k-m}\big)\big)\big)\big)
\label{eq:enum}
\end{equation}
as the generating function of $\ell$-slices with $2\leq \ell\leq k$ with a weight $g$ per inner face and a weight $\alpha^{\mathcal{L}(\ell-m)}$
whenever $\ell>  m$.

Assuming now $k\geq 3$ and $d$ in the range $2 \leq d\leq k-1$,  the generating function of $k$-slices with a weight $g$ per inner face and a weight $\alpha^{\mathcal{L}(d)}$ is given by
\begin{equation}
\underbrace{ \mathcal{K}^{(q)}\big( \mathcal{K}^{(q)}\big(\cdots \big( \mathcal{K}^{(q)}\big(}_{k-d\ \hbox{\scriptsize times}}\alpha^2\, T^{(q)}_{d}\big)\big)\big)\big)- \underbrace{ \mathcal{K}^{(q)}\big( \mathcal{K}^{(q)}\big(\cdots \big( \mathcal{K}^{(q)}\big(}_{k-d\ \hbox{\scriptsize times}}\alpha^2\, T^{(q)}_{d-1}\big)\big)\big)\big)\ .
\label{eq:enumbis}
\end{equation}
Indeed, the first term corresponds to $m=k-d$ in \eqref{eq:enum} hence enumerates $\ell$-slices with $2\leq \ell\leq k$  with a weight $\alpha^{\mathcal{L}(d+\ell-k)}$
if $d+\ell-k>0$ while the second term corresponds to taking $k\to k-1$ and $m=k-d$ in \eqref{eq:enum} hence enumerates $\ell$-slices with 
$2\leq \ell\leq k-1$  with \emph{the same weight} $\alpha^{\mathcal{L}(d+\ell-k)}$ if $d+\ell-k>0$. Taking the difference selects 
precisely $\ell$-slices with $\ell=k$, enumerated with a weight $\alpha^{\mathcal{L}(d)}$
(the condition $d>0$ is automatic since we assumed $d\geq 2$). 

Each of the term in the equation above may be computed as
follows: define $\lambda^{(q)}(\alpha;d)$ as the solution of the equation
\begin{equation*}
T_d^{(q)}(\lambda^{(q)}(\alpha;d))=\alpha^2\, T_d^{(q)}
\end{equation*}
or equivalently
\begin{equation}
\frac{(1-\lambda^{(q)}(\alpha;d)\, x^{d-1})(1-\lambda^{(q)}(\alpha;d)\, x^{d+4})}{(1-\lambda^{(q)}(\alpha;d)\, x^{d+1})(1-\lambda^{(q)}(\alpha;d)\,  x^{d+2})}=\alpha^2\, 
\frac{(1- x^{d-1})(1- x^{d+4})}{(1- x^{d+1})(1-  x^{d+2})}\ .
\label{eq:lambdaddef}
\end{equation}
This equation is quadratic in $\lambda^{(q)}(\alpha;d)$ and we have to pick the branch of solution satisfying $\lambda^{(q)}(1;d)=1$ \footnote{As already
noted, to be able to use property \eqref{eq:reclambda}, we have to ensure that $\frac{(1-\lambda^2\, x^{2 m+1})}{\left(1-\lambda\, x^{m}  \right) \left(1-\lambda\, x^{m+1}\right)}$ remains non-negative for $\lambda=\lambda^{(q)}(\alpha;d)$ and for all $m\leq d+1$. Since $x$ is in the range $0\leq x\leq 1$,
the most constraining requirement is for $m=d+1$, i.e.\ that $\frac{(1-\lambda^2\, x^{2 d+3})}{\left(1-\lambda\, x^{d+1}  \right) \left(1-\lambda\, x^{d+2}\right)}\geq 0$. Now this quantity changes sign upon changing $\lambda\to x^{-2d-3}/\lambda$, which precisely corresponds, in the equation 
for $\lambda$, to going from one branch of solution to the other. So only one of the solutions can be used, which by continuity is the branch
satisfying $\lambda^{(q)}(1;d)=1$ (the other branch satisfying $\lambda^{(q)}(1;d)=1/x^{2d+3}$ is not acceptable).}. This defines a \emph{unique}
value $\lambda^{(q)}(\alpha;d)$ and, from property \eqref{eq:reclambda}, we have
\begin{equation*}
\begin{split}
& \hskip -1.2cm \underbrace{ \big(\mathcal{K}^{(q)}\big( \mathcal{K}^{(q)}\big(\cdots \big( \mathcal{K}^{(q)}\big(}_{k-d\ \hbox{\scriptsize times}}\alpha^2\, T^{(q)}_{d}\big)\big)\big)\big)\big)= T^{(q)}_k(\lambda^{(q)}(\alpha;d))\\
& \hskip 4.cm =T^{(q)}\frac{(1-\lambda^{(q)}(\alpha;d)\, x^{k-1})(1-\lambda^{(q)}(\alpha;d)\, x^{k+4})}{(1-\lambda^{(q)}(\alpha;d)\, x^{k+1})(1-\lambda^{(q)}(\alpha;d)\,  x^{k+2})}\ .\\
\end{split}
\end{equation*}
By the same argument, we may compute the subtracted term in \eqref{eq:enumbis} and we find eventually that the generating function of $k$-slices with 
a weight $\alpha^{\mathcal{L}(d)}$ is given by
\begin{equation}
\begin{split}
&\hskip -1.cm Z^{(q)}(\alpha;d,k)= T^{(q)}_k(\lambda^{(q)}(\alpha;d))-T^{(q)}_{k-1}(\lambda^{(q)}(\alpha;d-1)) \\
& \quad \quad \ =
\frac{(1-x)  (1-x^2) \left(1+4x+x^2\right)}{\left(1+x+x^2\right)}  \times \\
&\hskip -1.cm  \times
  \frac{ x^{k-1}  (\lambda^{(q)}(\alpha;d-1)- \lambda^{(q)}(\alpha;d)\, x) \left(1- \lambda^{(q)}(\alpha;d-1) \lambda^{(q)}(\alpha;d)\, x^{2 k+2}\right)}{
  \left(1-\lambda^{(q)}(\alpha;d-1)\, x^k \right) \left(1- \lambda^{(q)}(\alpha;d-1)\, x^{k+1}\right) \left(1- \lambda^{(q)}(\alpha;d)\, x^{k+1}\right) \left(1- \lambda^{(q)}(\alpha;d)\, x^{k+2}\right)}\\
 \end{split}
 \label{eq:result}
\end{equation}
with  $\lambda^{(q)}(\alpha;d)$ as in \eqref{eq:lambdaddef}. This quantity is also the generating function of $k$-pointed-rooted quadrangulations
as we defined them,
with a fixed distance $d(v_0,v_1)=k\geq 3$, with a weight $g$ per face and a weight $\alpha^{\mathcal{L}(d)}$ where $\mathcal{L}(d)$
is the hull perimeter at some fixed distance $d$ ($2\leq d<k$) from $v_0$.
The first terms of the expansion in $g$ of $Z^{(q)}(\alpha;d,k)$ for the first allowed values of $k$ and $d$ are listed in Appendix B. 

This calculation is trivially generalized to control simultaneously the perimeters at two distances $d_1$ and $d_2$  with
$2\leq d_1\leq d_2 <k$.  We simply have to properly ``insert" the weight $\alpha_1$ at the $d_1$-th step of the recursion, then the weight
$\alpha_2$ at the $d_2$-th step. By doing so, we find that the generating function of $k$-pointed-rooted quadrangulations
(with a fixed distance $d(v_0,v_1)=k\geq 3$) with a weight $g$ per face and a weight $\alpha_1^{\mathcal{L}(d_1)}\alpha_2^{\mathcal{L}(d_2)}$ ($\mathcal{L}(d_1)$ and $\mathcal{L}(d_2)$
being the hull perimeters at respective distances $d_1$ and $d_2 $ from $v_0$) is, for $d_1\leq d_2$:
\begin{equation}
Z^{(q)}(\alpha_1,\alpha_2;d_1,d_2)\equiv T^{(q)}_k(\lambda^{(q)}(\alpha_1,\alpha_2;d_1,d_2))-T^{(q)}_{k-1}(\lambda^{(q)}(\alpha_1,\alpha_2;d_1-1,d_2-1))\ ,
 \label{eq:resultd1d2}
\end{equation}
with $T^{(q)}_k(\lambda)$ as in \eqref{eq:Tlambdak}, and where $\lambda^{(q)}(\alpha_1,\alpha_2;d_1,d_2)$ is defined as the solution
of 
\begin{equation}
\begin{split}
&\frac{(1-\lambda^{(q)}(\alpha_1,\alpha_2;d_1,d_2)\, x^{d_2-1})(1-\lambda^{(q)}(\alpha_1,\alpha_2;d_1,d_2)\, x^{d_2+4})}{(1-\lambda^{(q)}(\alpha_1,\alpha_2;d_1,d_2)\, x^{d_2+1})(1-\lambda^{(q)}(\alpha_1,\alpha_2;d_1,d_2)\,  x^{d_2+2})}\\
&\qquad \qquad \qquad \qquad  \qquad \qquad  \qquad =\alpha_2^2\, 
\frac{(1- \lambda^{(q)}(\alpha_1;d_1)\, x^{d_2-1})(1- \lambda^{(q)}(\alpha_1;d_1)\, x^{d_2+4})}{(1-\lambda^{(q)}(\alpha_1;d_1)\,  x^{d_2+1})(1- \lambda^{(q)}(\alpha_1;d_1)\,  x^{d_2+2})}\\
\end{split}
\label{eq:lambdad1d2def}
\end{equation}
with $\lambda^{(q)}(\alpha_1;d_1)$ defined as in \eqref{eq:lambdaddef}. Again the equation is quadratic in $\lambda^{(q)}(\alpha_1,\alpha_2;d_1,d_2)$
and we pick the branch of solution satisfying $\lambda^{(q)}(\alpha_1,1;d_1,d_2)=\lambda^{(q)}(\alpha_1;d_1)$.

\subsection{The case of triangulations}
\label{sec:triang}
\subsubsection{The slice recursion}
\begin{figure}
\begin{center}
\includegraphics[width=9cm]{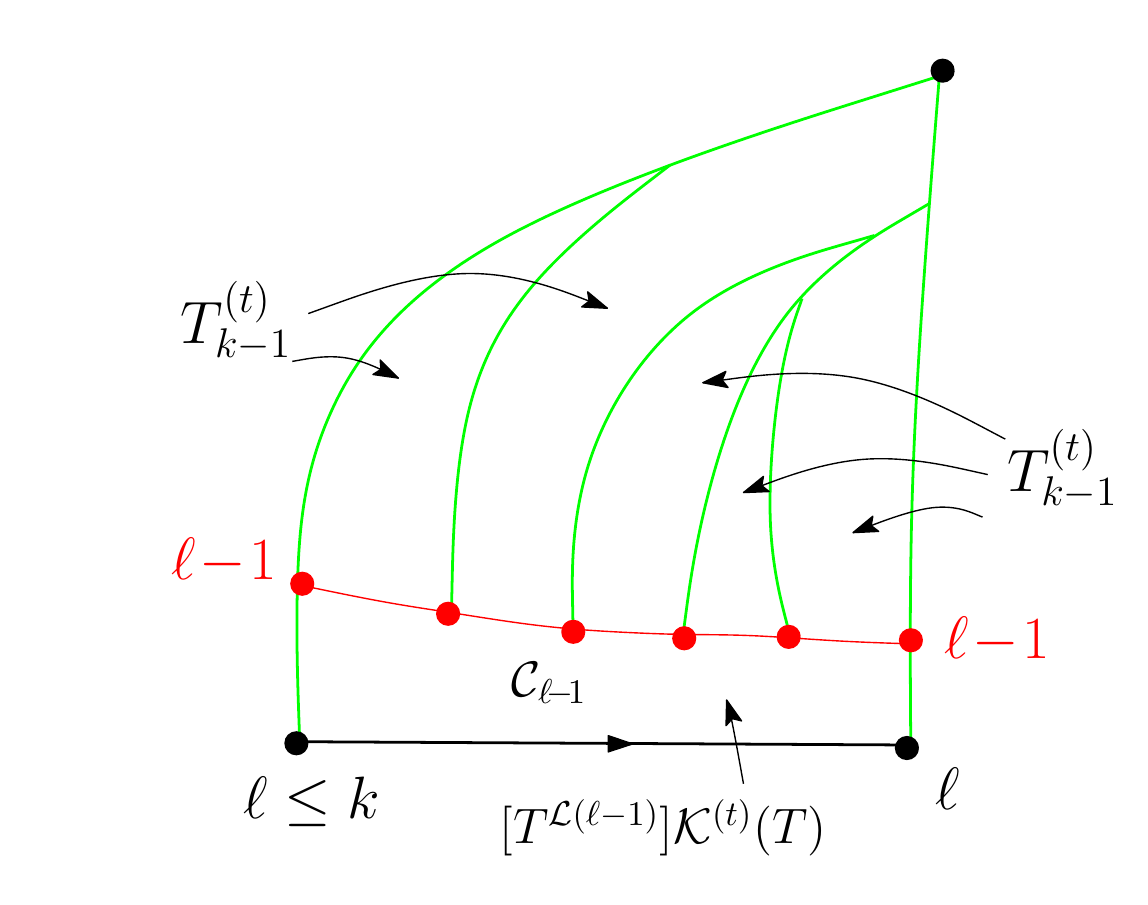}
\end{center}
\caption{Schematic representation of an $\ell$-isoslice and its decomposition by cutting along the hull boundary at distance $\ell-1$ from its apex. The domain 
$\mathcal{C}_{\ell-1}$ is enumerated by $[T^{\mathcal{L}}]\mathcal{K}^{(t)}(T)$ if the value of the hull perimeter $\mathcal{L}(\ell-1)$ is 
$\mathcal{L}$ while the complementary
domain is formed of $\mathcal{L}$ isoslices, each enumerated by $T^{(t)}_{k-1}$ (if we wish, for fixed $\mathcal{L}$, to count all the $\ell$-isoslices satisfying $1\leq \ell \leq k$). Here $\mathcal{L}=5$.}
\label{fig:constrhull3}
\end{figure}

Let us now discuss $k$-slices, as defined in Section~\ref{sec:definitions} in correspondence with $k$-pointed-rooted triangulations.
Again we consider the larger set of $\ell$-slices with $1\leq \ell \leq k$ and denote by $R^{(t)}_k\equiv R^{(t)}_k(g)$, $k\geq 1$ their 
generating function with a weight $g$ per inner triangle. 
The function $R^{(t)}_k$ is also the generating function of pointed-rooted triangulations with $1\leq d(v_0,v_1)\leq k$ and a weight $g$ per 
triangle\footnote{By 
pointed-rooted triangulations, we mean in general $\ell$-pointed-rooted triangulations with arbitrary $\ell\geq 1$. In particular,
the endpoint of the marked edge $e_1$ is at distance $1$ less from the origin than its first extremity.}.
The explicit expression for $R^{(t)}_k$ reads \cite{G15a}:
\begin{equation*}
R^{(t)}_k=R^{(t)} \frac{(1-x^{k})(1-x^{k+2})}{(1-x^{k+1})^2}\ , \qquad R^{(t)}=\frac{\sqrt{1+10x+x^2}}{1+x}\ .
\end{equation*}
where $x$ parametrizes $g$ through
\begin{equation*}
g=\frac{\sqrt{x(1+x)}}{(1+10x+x^2)^{3/4}}\ .
\end{equation*}
Again we fix $x$ univocally by imposing the extra condition $0\leq x\leq 1$: the generating functions are now well-defined for 
real $g$ in the range $0\leq g\leq 1/(2\cdot 3^{3/4})$.
The recursion relation involves, in addition to $R^{(t)}_k$, the generating function $T^{(t)}_k$ of \emph{$\ell$-isoslices} with $1\leq\ell\leq k$ 
and a weight $g$ per inner  triangle. The $\ell$-isoslices 
are defined exactly as $\ell$-slices except that both their right and left boundaries have the same length $\ell$
(in other fords, both extremities of the base are at distance $\ell$ from the apex -- see figure \ref{fig:constrhull3}).
Note that, as opposed to $\ell$-slices, $\ell$-isoslices are not related bijectively to some particular set of triangulations.
We have the explicit expression \cite{G15a}:
\begin{equation*}
T^{(t)}_k=T^{(t)} \frac{(1-x^{k})(1-x^{k+3})}{(1-x^{k+1})(1-x^{k+2})}\ , \qquad T^{(t)}=\sqrt{x}\, \frac{(1+10x+x^2)^{1/4}}{(1+x)^{3/2}}\ .
\end{equation*}

The recursion relation of \cite{G15a}, which fixes $R^{(t)}_k$ and $T^{(t)}_k$, may now be written as:
\begin{equation}
\begin{split}
& R^{(t)}_k=\mathcal{N}^{(t)}(T^{(t)}_{k-1})\,  \qquad T^{(t)}_k=\mathcal{K}^{(t)}(T^{(t)}_{k-1})\ ,\\
& \mathcal{N}^{(t)}(T)=\frac{R^{(t)}_1}{1-R^{(t)}_1\, T\, \Phi^{(t)}(T)}\ , \\
& \mathcal{K}^{(t)}(T)=\frac{(R^{(t)}_1)^2}{1-R^{(q)}_1\, T\, \Phi^{(t)}(T)}\ , \qquad  \Phi^{(t)}(T)\equiv \Phi^{(t)}(T,g)=\sum_{i\geq 3} h^{(t)}_{i}(g)\, T^{i-3}\ ,\\
\end{split}
\label{eq:recrelt}
\end{equation}
for $k\geq 1$ with $T^{(t)}_{0}=0$.
Here the quantities $h^{(t)}_{i}(g)$ denote appropriate generating functions defined in \cite{G15a} and whose explicit 
expression is not needed. Again we note that the kernels $\mathcal{N}^{(t)}(T)$ and $\mathcal{K}^{(t)}(T)$ are independent of $k$ 
while 
$R^{(t)}_k$ and $T^{(t)}_k$ depend on $k$ only through the variable $x^k$. 
A before, we immediately deduce that, if we set
\begin{equation}
R^{(t)}_k(\lambda)=R^{(t)} \frac{(1-\lambda\, x^{k})(1-\lambda\, x^{k+2})}{(1-\lambda\, x^{k+1})^2}\ , \qquad
T^{(t)}_k(\lambda) =T^{(t)} \frac{(1-\lambda\, x^{k})(1-\lambda\, x^{k+3})}{(1-\lambda\, x^{k+1})(1-\lambda\, x^{k+2})} \ , 
\label{eq:RTlambdak}
\end{equation} 
then 
\begin{equation}
R^{(t)}_k(\lambda)=\mathcal{N}^{(t)}(T^{(t)}_{k-1}(\lambda))\ , \qquad 
T^{(t)}_k(\lambda)=\mathcal{K}^{(t)}(T^{(t)}_{k-1}(\lambda))
\label{eq:reclambdat}
\end{equation}
(for $\lambda$ close enough to $1$).

\subsubsection{Controlling the hull perimeter}
Again the recursion relation is intimately linked to the notion of hull perimeter. Let us recall how it works for triangulations. 
We consider an $\ell$-isoslice with left boundary length $\ell$ between $1$ and $k$ (as enumerated by $T^{(t)}_k$)
and cut it along its so-called dividing line (as defined in \cite{G15a}) which is nothing but
the hull boundary at distance $\ell-1$ from $v_0$ in the $\ell$-isoslice, defined exactly as in Section \ref{sec:definitions}
(see figure \ref{fig:constrhull3} for an illustration)\footnote{When $\ell=1$, the hull boundary at distance $\ell-1=0$ must be understood as reduced to the single
vertex $v_0$, having length $0$.}. Thanks to this cutting, we deduce that the generating
function $T^{(t)}_k$ is the product of the generating function of the domain $\mathcal{C}_{\ell-1}$ times that of the domain $\mathcal{H}_{\ell-1}$.
For a fixed value $\mathcal{L}$ of the hull perimeter $\mathcal{L}(\ell-1)$, the first generating function is simply:
\begin{equation*}
[T^{\mathcal{L}}]\mathcal{K}^{(t)}(T)
\end{equation*}
while the second generating function reads:
\begin{equation*}
\left(T^{(t)}_{k-1}\right)^{\mathcal{L}}\ .
\end{equation*}
This is because the domain $\mathcal{H}_{\ell-1}$ may be decomposed into exactly $\mathcal{L}$ sub-isoslices with left boundary lengths $\ell_i$, 
$1\leq i\leq \mathcal{L}$, satisfying $1\leq \ell_i\leq \ell-1$, hence $1\leq \ell_i\leq k-1$ when considering all possible values of $\ell$. As before, 
these sub-isoslices are obtained by cutting along the leftmost shortest paths to $v_0$
from the ${\mathcal{L}}-1$ internal vertices of the hull boundary at distance $\ell-1$. Each of the ${\mathcal{L}}$ sub-isoslices contributes
a factor $T^{(t)}_{k-1}$ to the generating function.
Summing over all values of $\mathcal{L}$ yields the desired recursion relation \eqref{eq:recrelt}. To assign a weight $\alpha^{\mathcal{L}}$ 
to our $\ell$-isoslices, we simply need to replace $T^{(t)}_{k-1}$ by 
$\alpha\, T^{(t)}_{k-1}$ at the step $k-1\to k$ of the recursion so that the generating function of $\ell$-isoslices with $1\leq \ell\leq k$, with a weight $g$ per inner face and a weight 
$\alpha^{\mathcal{L}(\ell-1)}$ is simply
\begin{equation*}
 \mathcal{K}^{(t)}(\alpha\, T^{(t)}_{k-1})
\end{equation*}
(again, $\ell$-isoslices with $\ell=1$ must be understood as having $\mathcal{L}(\ell-1)=0$, hence are enumerated with a weight $\alpha^{\mathcal{L}(\ell-1)}=1$).

\begin{figure}
\begin{center}
\includegraphics[width=8cm]{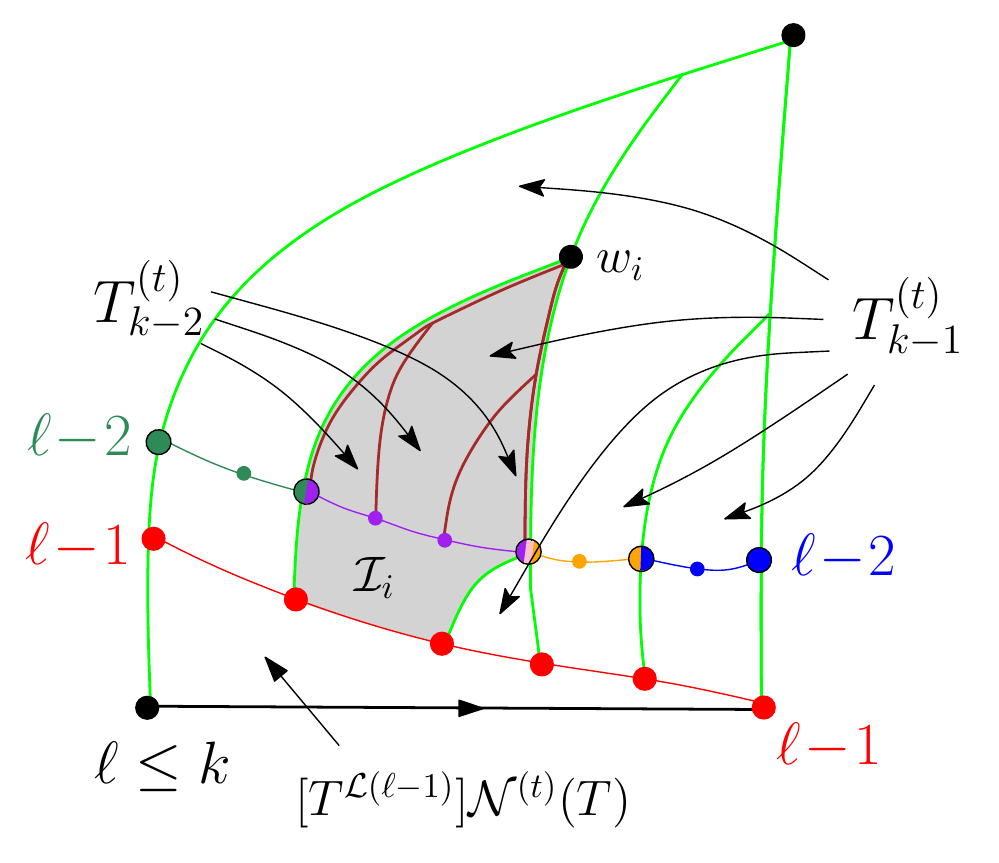}
\end{center}
\caption{Schematic representation of an $\ell$-slice and its decomposition by cutting along the hull boundary at distance $\ell-1$ from its apex. The domain 
below the hull boundary is now enumerated by $[T^{\mathcal{L}}]\mathcal{N}^{(t)}(T)$ if the value of the hull perimeter $\mathcal{L}(\ell-1)$ is $\mathcal{L}$ while the complementary
domain is formed of $\mathcal{L}$ isoslices $\mathcal{I}_i$ of left boundary length $\ell_i$, exactly as in figure \ref{fig:constrhull3}, each enumerated by $T^{(t)}_{k-1}$ (if we wish, for fixed $\mathcal{L}$, to enumerate all $\ell$-slices with $1\leq \ell \leq k$). The hull boundary at distance $\ell-2$ from the apex is
obtained by concatenating the hull boundaries at distance $\ell_i-1$ from their apex $w_i$ of all the sub-isoslices $\mathcal{I}_i$.}
\label{fig:constrhull4}
\end{figure} 
\vskip .3cm
As for the relation giving $R^{(t)}_{k}$ in \eqref{eq:recrelt}, it follows from a similar decomposition of $\ell$-slices with left boundary length $\ell$ between $1$ and $k$. The hull $\mathcal{H}_{\ell-1}$ is characterized by exactly the same distance constraints as for isoslices, 
hence yields the same generating function $\left(T^{(t)}_{k-1}\right)^{\mathcal{L}}$ if $\mathcal{L}(\ell-1)$ has a fixed value ${\mathcal{L}}$.
Only the domain $\mathcal{C}_{\ell-1}$ is modified and has then generating function
$[T^{\mathcal{L}}]\mathcal{N}^{(t)}(T)$ (see figure \ref{fig:constrhull4}). This leads to the desired relation in \eqref{eq:recrelt}. Again, we may easily assign a weight $\alpha^{\mathcal{L}}$ by multiplying $T^{(t)}_{k-1}$ by $\alpha$. In other words, the generating function of $\ell$-slices with $1\leq \ell\leq k$, with a weight $g$ per inner face and a weight 
$\alpha^{\mathcal{L}(\ell-1)}$ is
\begin{equation*}
 \mathcal{N}^{(t)}(\alpha\, T^{(t)}_{k-1})
\end{equation*}
($\ell$-slices with $\ell=1$ are enumerated with a weight $\alpha^{\mathcal{L}(\ell-1)}=1$).

Repeating the argument of  Section~\ref{sec:quad}, we obtain, for $1\leq m<k$,  
the generating functions of, respectively,  $\ell$-slices and $\ell$-isoslices with $1\leq \ell\leq k$ with a weight $g$ per inner face and a weight $\alpha^{\mathcal{L}(\ell-m)}$ whenever $\ell \geq m$, namely:
\begin{equation}
\begin{split}
&\mathcal{N}^{(t)} \big(\underbrace{\mathcal{K}^{(t)}\big(\cdots \big( \mathcal{K}^{(t)}\big(}_{m-1\ \hbox{\scriptsize times}}\alpha\, T^{(t)}_{k-m}\big)\big)\big)\big) \qquad \hbox{and} \\
&\underbrace{ \mathcal{K}^{(t)}\big( \mathcal{K}^{(t)}\big(\cdots \big( \mathcal{K}^{(t)}\big(}_{m\ \hbox{\scriptsize times}}\alpha\, T^{(t)}_{k-m}\big)\big)\big)\big)  \qquad \hbox{respectively}\\
\end{split}
\label{eq:enumt}
\end{equation}
(see figure \ref{fig:constrhull4} for an illustration of the first identity when $m=2$).

If we now take $k\geq 2$ and $d$ in the range $1 \leq d\leq k-1$, the generating function of $k$-slices with a weight $g$ per inner face and a weight $\alpha^{\mathcal{L}(d)}$ is 
\begin{equation*}
\mathcal{N}^{(t)}\big(\underbrace{\mathcal{K}^{(t)}\big(\cdots \big( \mathcal{K}^{(t)}\big(}_{k-d-1\ \hbox{\scriptsize times}}\alpha\, T^{(t)}_{d}\big)\big)\big)\big)- \mathcal{N}^{(t)} \big(\underbrace{\mathcal{K}^{(t)}\big(\cdots \big( \mathcal{K}^{(t)}\big(}_{k-d-1 \hbox{\scriptsize times}}\alpha\, T^{(t)}_{d-1}\big)\big)\big)\big)\ .
\label{eq:enumtbis}
\end{equation*}
As for quadrangulations, both terms in the equation may be computed upon defining $\lambda^{(t)}(\alpha;d)$ as the solution of the equation
\begin{equation*}
T_d^{(t)}(\lambda^{(t)}(\alpha;d))=\alpha\, T_d^{(t)}\ ,
\end{equation*}
namely
\begin{equation}
\frac{(1-\lambda^{(t)}(\alpha;d)\, x^{d})(1-\lambda^{(t)}(\alpha;d)\, x^{d+3})}{(1-\lambda^{(t)}(\alpha;d)\, x^{d+1})(1-\lambda^{(t)}(\alpha;d)\,  x^{d+2})}=\alpha\, 
\frac{(1- x^{d})(1- x^{d+3})}{(1- x^{d+1})(1-  x^{d+2})}
\label{eq:lambdaddeft}
\end{equation}
where we pick the branch of solution satisfying $\lambda^{(t)}(1;d)=1$. From property \eqref{eq:reclambdat}, we have
\begin{equation*}
\begin{split}
& \hskip -1.2cm \mathcal{N}^{(t)}\big( \underbrace{\mathcal{K}^{(t)}\big(\cdots \big( \mathcal{K}^{(t)}\big(}_{k-d-1\ \hbox{\scriptsize times}}\alpha\, T^{(q)}_{d}\big)\big)\big)\big)= R^{(t)}_k(\lambda^{(t)}(\alpha;d))\\
& \hskip 3.6cm =R^{(t)}\frac{(1-\lambda^{(t)}(\alpha;d)\, x^{k})(1-\lambda^{(t)}(\alpha;d)\, x^{k+2})}{(1-\lambda^{(t)}(\alpha;d)\, x^{k+1})^2}\ .\\
\end{split}
\end{equation*}
so that the generating function of $k$-slices with a weight $\alpha^{\mathcal{L}(d)}$ eventualy reads
\begin{equation}
\begin{split}
& \hskip -1.cm Z^{(t)}(\alpha;d,k)= R^{(t)}_k(\lambda^{(t)}(\alpha;d))-R^{(t)}_{k-1}(\lambda^{(t)}(\alpha;d-1)) \\ & \qquad =
\frac{(1-x)^2  \sqrt{1+10x+x^2}}{(1+x)}  \times \\
& \qquad \quad  \times
 \frac{ x^{k-1}  (\lambda^{(t)}(\alpha;d-1)- \lambda^{(t)}(\alpha;d)\, x) \left(1- \lambda^{(t)}(\alpha;d-1) \lambda^{(t)}(\alpha;d)\, x^{2 k+1}\right)}{
 \left(1-\lambda^{(t)}(\alpha;d-1)\, x^k \right)^2 \left(1- \lambda^{(t)}(\alpha;d)\, x^{k+1}\right)^2 }\\
 \end{split}
 \label{eq:resultt}
\end{equation}
with $\lambda^{(t)}(\alpha;d)$ as in \eqref{eq:lambdaddeft}. This is also the generating function of $k$-pointed-rooted triangulations as we defined them,
with a fixed distance $d(v_0,v_1)=k\geq 2$, with a weight $g$ per triangle and a weight $\alpha^{\mathcal{L}(d)}$ where $\mathcal{L}(d)$
is the hull perimeter at some fixed distance $d$ ($1\leq d<k$) from $v_0$.
The first terms of the expansion in $g$ of $Z^{(t)}(\alpha;d,k)$ for the first allowed values of $k$ and $d$ are listed in Appendix B. 

Again we may impose a simultaneous control on the perimeters at two distances $d_1$ and $d_2$  ($1\leq d_1\leq d_2 <k$).  
The generating function of $k$-pointed-rooted triangulations (with a fixed distance $d(v_0,v_1)=k\geq 2$) with a weight $g$ per face and a weight $\alpha_1^{\mathcal{L}(d_1)}\alpha_2^{\mathcal{L}(d_2)}$ ($\mathcal{L}(d_1)$ and $\mathcal{L}(d_2)$
being the hull perimeters at respective distances $d_1$ and $d_2 $ from $v_0$) is, for $d_1\leq d_2$:
\begin{equation}
Z^{(t)}(\alpha_1,\alpha_2;d_1,d_2)\equiv R^{(t)}_k(\lambda^{(t)}(\alpha_1,\alpha_2;d_1,d_2))-R^{(t)}_{k-1}(\lambda^{(t)}(\alpha_1,\alpha_2;d_1-1,d_2-1))\ ,
 \label{eq:resultd1d2t}
\end{equation}
with $R^{(t)}_k(\lambda)$ as in \eqref{eq:RTlambdak}, and where $\lambda^{(t)}(\alpha_1,\alpha_2;d_1,d_2)$ is defined as the solution
of 
\begin{equation}
\begin{split}
&\frac{(1-\lambda^{(t)}(\alpha_1,\alpha_2;d_1,d_2)\, x^{d_2})(1-\lambda^{(t)}(\alpha_1,\alpha_2;d_1,d_2)\, x^{d_2+3})}{(1-\lambda^{(t)}(\alpha_1,\alpha_2;d_1,d_2)\, x^{d_2+1})(1-\lambda^{(t)}(\alpha_1,\alpha_2;d_1,d_2)\,  x^{d_2+2})}\\
&\qquad \qquad \qquad \qquad  \qquad \qquad  \qquad =\alpha_2\, 
\frac{(1- \lambda^{(t)}(\alpha_1;d_1)\, x^{d_2})(1- \lambda^{(t)}(\alpha_1;d_1)\, x^{d_2+3})}{(1-\lambda^{(t)}(\alpha_1;d_1)\,  x^{d_2+1})(1- \lambda^{(t)}(\alpha_1;d_1)\,  x^{d_2+2})}\\
\end{split}
\label{eq:lambdad1d2def}
\end{equation}
with $\lambda^{(t)}(\alpha_1;d_1)$ defined as in \eqref{eq:lambdaddeft}. The correct branch of solutions for $\lambda^{(t)}(\alpha_1,\alpha_2;d_1,d_2)$
is that satisfying $\lambda^{(t)}(\alpha_1,1;d_1,d_2)=\lambda^{(t)}(\alpha_1;d_1)$.

\section{The limit of large maps}
\label{sec:largemaps}
\subsection{Singularity analysis}
\label{sec:singularity}
We now have at our disposal all the required generating functions. Let us recall how to extract from these quantities
the desired expectation values and probability distributions. Our results of Section~\ref{sec:mainresults} apply to the ensemble of uniformly 
drawn $k$-pointed-rooted quadrangulations (respectively triangulations) having a fixed number  $N$ of faces and a fixed value $k$ of the distance between their origin $v_0$ and the first extremity $v_1$ of their marked edge $e_1$. In our generating functions, $k$ is already fixed
but, in order to have a fixed $N$, we must in principle extract the coefficient of $g^N$ in their expansion in $g$. In Section~\ref{sec:mainresults},
we specialized our results to 
the limit $N\to \infty$: the asymptotic behavior of the coefficient of $g^N$ is then directly encoded in the singular behavior 
of the generating functions. In the case of quadrangulations, a singularity occurs when $g$ approaches its ``critical value" $1/12$,
corresponding to the limit $x\to 1$.
For triangulations,  the singularity is for $g\to 1/(2\cdot 3^{3/4})$. The singular behavior is best captured by setting
\begin{equation*}
\begin{split}
&g=\frac{1}{12}\left(1-\epsilon^4\right)\qquad \hbox{for quadrangulations}\ ,\\
&g=\frac{1}{2\cdot 3^{3/4}}\left(1-\epsilon^4\right)\qquad \hbox{for triangulations}\ ,\\
\end{split}
\end{equation*}
 and expanding our various generating functions around $\epsilon=0$. 

Let us now discuss in details the case of quadrangulations. The small $\epsilon$ expansion of $Z^{(q)}(\alpha;d,k)$ 
takes the form: 
\begin{equation}
Z^{(q)}(\alpha;d,k)=A_0^{(q)}(\alpha;d,k)+A_2^{(q)}(\alpha;d,k)\epsilon^4+A_3^{(q)}(\alpha;d,k)\epsilon^6+O(\epsilon^8)\ ,
\label{eq:Zexpan}
\end{equation}
where we note in particular \emph{the absence of a term of order $\epsilon^2$}. The constant term and the term of 
order $\epsilon^4=1-12\, g$ being regular, the most singular part of  $Z^{(q)}(\alpha;d,k)$ is therefore
\begin{equation*}
Z^{(q)}(\alpha;d,k)\vert_{\rm sing.}= A_3^{(q)}(\alpha;d,k)\ (1-12\, g)^{3/2}
\end{equation*}
and, by a standard result, we deduce the large $N$ behavior
\begin{equation*}
[g^N]Z^{(q)}(\alpha;d,k)\sim \frac{3}{4}\  \frac{12^N}{N^{5/2}}\times A_3^{(q)}(\alpha;d,k)\ .
\end{equation*}

This gives the large $N$ asymptotics of the ``reduced" generating function (with $\alpha$ as only left variable) 
of $k$-pointed-rooted quadrangulations with a fixed number $N$ of faces,
a fixed distance $d(v_0,v_1)=k$ and a weight $\alpha^{\mathcal{L}(d)}$.
To get the expectation value of $\alpha^{\mathcal{L}(d)}$
in the ensemble of $k$-pointed-rooted quadrangulations with fixed $N$ and $k$, we must divide this generating function 
by the cardinal of this ensemble. The latter is easily obtained by taking $\alpha=1$, in which case $Z^{(q)}(\alpha;d,k)$
and $A_3^{(q)}(\alpha;d,k)$ do not depend on $d$ (recall that $\lambda(1;d)=1$ independently of $d$).
The large $N$ behavior of the \emph{number} of $k$-pointed-rooted quadrangulations with fixed $N$ and $k$ is thus
\begin{equation*}
[g^N]Z^{(q)}(1;d,k)\sim \frac{3}{4}\  \frac{12^N}{N^{5/2}}\times A_3^{(q)}(1;d,k)
\end{equation*}
irrespectively of $d$. The expectation value of $\alpha^{\mathcal{L}(d)}$ in our ensemble is thus simply
\begin{equation}
W_k^{(q)}(\alpha;d)\equiv E_k(\alpha^{\mathcal{L}(d)})= \frac{A_3^{(q)}(\alpha;d,k)}{A_3^{(q)}(1;d,k)}\ .
\label{eq:Wkd}
\end{equation}
The computation of $W_k^{(q)}(\alpha;d)$, although straightforward, is quite cumbersome and we do not give its full expression here. 
As we discussed in Section \ref{sec:mainresults}, our main interest is the limit of large $k$. If we send $k\to \infty$, \emph{keeping $d$ finite},
the expression of $W_k^{(q)}(\alpha;d)$ drastically simplifies and we find
\begin{equation}
\begin{split}
&\hskip -1.cm W_\infty^{(q)}(\alpha;d)= \frac{1}{2} \left(
-\frac{\sqrt{(d-2)^2(d+3)^2 \alpha ^4-26 (d-2) d (d+1)
  (d+3) \alpha ^2+25 d^2 (d+1)^2}}{d (d+1) \left(1-\alpha ^2\right)+6 \alpha
   ^2}
   \right.\\
   & +\left. 
   \frac{\sqrt{(d-1)^2 (d+4)^2 \alpha ^4-26 (d-1) (d+1) (d+2) (d+4) \alpha
   ^2+25 (d+1)^2 (d+2)^2}}{(d+1) (d+2) \left(1-\alpha ^2\right)+6 \alpha
   ^2}
   \right)\ . \\
   \end{split}
 \label{eq:Wqinf}
\end{equation}
In particular, the average length $\mathcal{L}(d)$ is simply
\begin{equation*}
E_\infty (\mathcal{L}(d))=\frac{\partial}{\partial \alpha} W_\infty^{(q)}(\alpha;d)\vert_{\alpha=1}=\frac{2 (d+1)^2 \left(3 d^2+6 d-4\right)}{3 (2d+1)(2d+3)}\ .
\end{equation*}
For large $d$, this average length scales like $d^2/2= (3c/2)\times d^2$ with $c=1/3$. Introducing as in Section \ref{sec:mainresults}
the rescaled variable $L(d)\equiv \mathcal{L}(d)/d^2$,
we have
\begin{equation*}
\lim_{d\to \infty} E_\infty (e^{-\tau L(d)})= \lim_{d\to \infty} W_\infty^{(q)}(e^{-\frac{\tau}{d^2}};d)= \frac{1}{(1+c\, \tau)^{3/2}}\
\end{equation*}  
with $c=1/3$, as announced in eq.~\eqref{eq:CLG}.

Returning to a finite value of $k$, we give in Appendix C the general expression for the average length $\mathcal{L}(d)$, i.e. the expression
of 
\begin{equation*}
E_k(\mathcal{L}(d))=\frac{\partial}{\partial \alpha} W_k^{(q)}(\alpha;d)\vert_{\alpha=1}\ .
\end{equation*} 
From this expression, we deduce, by considering the limit of both $d$ and $k$ large, with a fixed value of the ratio $u=d/k$ ($0<u<1$), the average value $L_{\rm av}(u)$ of rescaled length $L(d)$ , namely:
\begin{equation*}
L_{\rm av}(u) \equiv \lim_{k\to \infty} \frac{1}{(k\, u)^2} E_k(\mathcal{L}(k\, u))= \frac{3c}{2}(1+u-3 u^6+u^7)
\end{equation*}
with $c=1/3$, as announced in Section \ref{sec:mainresults}.
\vskip .3cm
We may perform similar calculations for triangulations, taking as starting point the quantity $Z^{(t)}(\alpha;d,k)$ of eq.~\eqref{eq:resultt}.
Following the same lines as above, we obtain an expression for
\begin{equation*}
W_k^{(t)}(\alpha;d)\equiv E_k(\alpha^{\mathcal{L}(d)})
\end{equation*}
in the ensemble of $k$-pointed-rooted triangulations with fixed $k$ and $N$, in the asymptotic limit where $N\to \infty$. 
In particular, we now get
\begin{equation}
\begin{split}
&\hskip -1.cm  W_\infty ^{(t)}(\alpha;d)=\frac{1}{2} \left(
-\frac{\sqrt{(d-1) (d+2) ((d-1) (d+2) (9-\alpha) (1-\alpha )-20 \alpha +36)+36}}{(d-1) (d+2) (1-\alpha
   )+2}
   \right.\\
   & \qquad \qquad \qquad \qquad \qquad +\left. 
\frac{\sqrt{d (d+3) (d (d+3) (9-\alpha) (1-\alpha)-20 \alpha +36)+36}}{d (d+3) (1-\alpha)+2}   
   \right)\ ,\\
   \end{split}
    \label{eq:Wtintf}
\end{equation}
from which we extract
\begin{equation*}
E_\infty (\mathcal{L}(d))=\frac{3 d (d+1)^2 (d+2)}{(2d+1)(2d+3)}\ ,
\end{equation*}
which now scales at large $d$ as $3 d^2/4= (3c/2)\times d^2$ with $c=1/2$. The explicit formula of $E_k(\mathcal{L}(d))$ for
finite $k$ and $d$ is presented in Appendix C.
We also easily compute as above the values of $ \lim_{d\to \infty} E_\infty (e^{-\tau L(d)})$ and $L_{\rm av}(u)$ for triangulations, whose expressions are 
exactly the same as for quadrangulations, except for the value of $c$, now equal to $1/2$.

\subsection{A shortcut in the calculations}
\label{sec:shortcut}
The calculations above were performed for finite $d$ and $k$, leading to some non-universal quantities
with rather involved expressions from which we then extracted some of the universal laws at $d,k\to \infty$ announced
in Section~\ref{sec:mainresults}, such as the expressions
for $L_{\rm av}(u)$ and for $\lim_{d\to \infty} E_\infty (e^{-\tau L(d)})$. To obtain in a more systematic way 
the results of Section~\ref{sec:mainresults}, we may simplify our calculations and incorporate {\it ab initio} the fact 
that $k$ and $d$ are eventually supposed to be large. 

\subsubsection{A more direct computation of the probability density for $L(d)$ when $k\to \infty$}
We start again with the case of quadrangulations. All our calculations rely on the expression of $A_3^{(q)}(\alpha;d,k)$ in \eqref{eq:Zexpan}. The large $k$ behavior of this coefficient
may be deduced from the so-called \emph{scaling limit} where we let $\epsilon\to 0$ and $k\to \infty$ as 
\begin{equation*}
k=\frac{K}{\epsilon}
\end{equation*} 
with $K$ finite. More precisely, when $\epsilon\to 0$, we have the expansions
\begin{equation*}
\begin{split}
& x=1-\sqrt{6}\, \epsilon+O(\epsilon^2)\\
&\lambda^{(q)}(\alpha;d)=1-\Lambda^{(q)}(\alpha;d)\epsilon+O(\epsilon^2)\\
\end{split}
\end{equation*}
where $\Lambda^{(q)}(\alpha;d)$ is obtained by expanding \eqref{eq:lambdaddef} at leading order in $\epsilon$, namely as the solution of
\begin{equation}
\frac{\left(\Lambda^{(q)}(\alpha;d)\right)^2+\sqrt{6} (2 d+3) \Lambda^{(q)}(\alpha;d) +6 (d-1) (d+4)}{\left(\Lambda^{(q)}(\alpha;d)\right)^2+\sqrt{6} (2
   d+3) \Lambda^{(q)}(\alpha;d) +6 (d+1) (d+2)}=\alpha^2\, \frac{(d-1) (d+4)}{(d+1)( d+2)}
\label{eq:Lambdadef}
\end{equation}
which vanishes when $\alpha=1$. We then have the scaling behavior when $\epsilon\to 0$:
\begin{equation*}
\begin{split}
&Z^{(q)}\left(\alpha;d,\frac{K}{\epsilon}\right)\sim \mathcal {Z}^{(q)}(\alpha;d,K)\,  \epsilon^3 \\ 
&\mathcal {Z}^{(q)}(\alpha;d,K)=  \left(  \frac{24 e^{\sqrt{6} K} \left(1+e^{\sqrt{6} K}\right) \left( \Lambda^{(q)}(\alpha;d)-  \Lambda^{(q)}(\alpha;d-1)+\sqrt{6}\right)}{\left(e^{\sqrt{6} K}-1\right)^3}\right) \ .
\end{split} 
\end{equation*}
To be consistent with the expansion \eqref{eq:Zexpan}, the coefficient $A_i^{(q)}(\alpha;d,k)$ must behave at large $k$ as $k^{2i-3}$ 
and more precisely
\begin{equation*}
A_i^{(q)}(\alpha;d,k)\underset{k\to \infty}{\sim} k^{2i-3} \times [K^{2i-3}]\mathcal {Z}^{(q)}(\alpha;d,K)\ .
\end{equation*}
Expanding $\mathcal {Z}^{(q)}(\alpha;d,K)$ at small $K$ and extracting the term of order $K^3$, we thus deduce
\begin{equation*}
A_3^{(q)}(\alpha;d,k)\underset{k\to \infty}{\sim} k^{3} \times \frac{2}{7} \sqrt{\frac{2}{3}} \left( \Lambda^{(q)}(\alpha;d)-  \Lambda^{(q)}(\alpha;d-1)+\sqrt{6}\right)
\end{equation*}
and, from \eqref{eq:Wkd},
\begin{equation}
E_k(\alpha^{\mathcal{L}(d)})=W_k^{(q)}(\alpha;d)
\underset{k\to \infty}{\sim} \frac{\Lambda^{(q)}(\alpha;d)-  \Lambda^{(q)}(\alpha;d-1)+\sqrt{6}}{\sqrt{6}}
\label{eq:avkinf}
\end{equation}
from which  \eqref{eq:Wqinf} follows immediately. When $d\to \infty$, we easily get from \eqref{eq:Lambdadef}
that 
\begin{equation*}
\begin{split}
& \hskip -1.cm \Lambda^{(q)}(e^{-\frac{\tau}{d^2}};d)=\sqrt{6}\, d \left(\frac{1}{\sqrt{1+c\, \tau}}-1\right)+3 \sqrt{\frac{3}{2}}
   \left(\frac{1}{\left(1+c\, \tau\right)^{3/2}}-1\right)+O\left(\frac{1}{d}\right) \\
&  \hskip -1.cm \Lambda^{(q)}(e^{-\frac{\tau}{d^2}};d-1)= \sqrt{6}\, d \left(\frac{1}{\sqrt{1+c\, \tau}}-1\right)+\sqrt{\frac{3}{2}}
   \left(\frac{1}{\left(1+c\, \tau\right)^{3/2}}-1\right)+O\left(\frac{1}{d}\right) \\
   \end{split}
\end{equation*}
with $c=1/3$ and \eqref{eq:avkinf} immediately implies \eqref{eq:CLG} and, by some inverse Laplace transform, the expression \eqref{eq:Pinf} of 
the probability density for $L(d)$ when $d$ is large. In conclusion, the above approach based on the scaling limit provides a more direct proof of our first main result, here for quadrangulations.

\subsubsection{Computation of the joint probability density for $L(d_1)$ and $L(d_2)$ when $k \to\infty$}
By a straightforward extension of the above analysis, we find that, for $d_1\leq d_2$, 
\begin{equation}
E_k(\alpha_1^{\mathcal{L}(d_1)}\alpha_2^{\mathcal{L}(d_2)})
\underset{k\to \infty}{\sim} \frac{\Lambda^{(q)}(\alpha_1,\alpha_2;d_1,d_2)-  \Lambda^{(q)}(\alpha_1,\alpha_2;d_1-1,d_2-1)+\sqrt{6}}{\sqrt{6}}
\label{eq:Ekalpha12}
\end{equation}
where $\Lambda^{(q)}(\alpha_1,\alpha_2;d_1,d_2)$ is fixed by 
\begin{equation*}
\begin{split}
& \hskip -1.cm \frac{\left(\Lambda^{(q)}(\alpha_1,\alpha_2;d_1,d_2)\right)^2+\sqrt{6} (2 d_2+3) \Lambda^{(q)}(\alpha_1,\alpha_2;d_1,d_2) +6 (d_2-1) (d_2+4)}{\left(\Lambda^{(q)}(\alpha_1,\alpha_2;d_1,d_2)\right)^2+\sqrt{6} (2
   d_2+3) \Lambda^{(q)}(\alpha_1,\alpha_2;d_1,d_2) +6 (d_2+1) (d_2+2)}= \\
   &\qquad \qquad \qquad  \alpha_2^2\, \frac{\left(\Lambda^{(q)}(\alpha_1;d_1)\right)^2+\sqrt{6} (2 d_2+3) \Lambda^{(q)}(\alpha_1;d_1) +6 (d_2-1) (d_2+4)}{\left(\Lambda^{(q)}(\alpha_1;d_1)\right)^2+\sqrt{6} (2
   d_2+3) \Lambda^{(q)}(\alpha_1;d_1) +6 (d_2+1) (d_2+2)}\\
   \end{split}
\end{equation*}
with $\Lambda^{(q)}(\alpha_1;d_1)$ defined by \eqref{eq:Lambdadef} (we pick the solution such that $\Lambda^{(q)}(\alpha_1,1;d_1,d_2)=
\Lambda^{(q)}(\alpha_1;d_1)$).  Setting $d_1=d$, $d_2=v\, d$ ($v\geq 1$) and letting $d\to \infty$, we then have
\begin{equation*}
\begin{split}
& \hskip -1.cm \Lambda^{(q)}(e^{-\frac{\tau_1}{d^2}},e^{-\frac{\tau_2}{(v\, d)^2}};d,v\, d)= \sqrt{6}\,  d\,  v\, \left( \frac{ 1+(v-1) \omega_1}{\sqrt{v^2\, \omega_1
   ^2+c \,\tau_2 (1+(v-1) \omega_1)^2}}-1\right)
\\ &\hskip 3.cm  +3 \sqrt{\frac{3}{2}} \left(\frac{v^3}{\left(v^2\, \omega_1^2+c\, \tau_2 (1+(v-1)
   \omega_1)^2\right)^{3/2}}-1\right)+O\left(\frac{1}{d}\right) 
 \\
&  \hskip -1.cm \Lambda^{(q)}(e^{-\frac{\tau_1}{d^2}},e^{-\frac{\tau_2}{(v\, d)^2}};d-1,v\, d-1)=\sqrt{6}\,  d\,  v\, \left( \frac{ 1+(v-1) \omega_1}{\sqrt{v^2\, \omega_1
   ^2+c \,\tau_2 (1+(v-1) \omega_1)^2}}-1\right)
\\ &\hskip 3.cm + \sqrt{\frac{3}{2}} \left(\frac{v^3}{\left(v^2\, \omega_1^2+c\, \tau_2 (1+(v-1)
   \omega_1)^2\right)^{3/2}}-1\right))+O\left(\frac{1}{d}\right) \\
   & \hskip -1.cm \hbox{where}\ \omega_1\equiv \sqrt{1+c\, \tau_1}\ ,\\
   \end{split}
\end{equation*}
with $c=1/3$.
Plugging these expansions in \eqref{eq:Ekalpha12}, we deduce the announced result \eqref{eq:zuniv}. 
Knowing \eqref{eq:zuniv}, it is not difficult to deduce the expression \eqref{eq:PL1L2} for the joint probability
density $\mathcal{P}(L_1,L_2;v)$ by performing a double inverse Laplace transform. The details of this transformation
are presented in Appendix A. This proves our third main result
for quadrangulations.

\subsubsection{Similar results for triangulations}
In the case of triangulations, we find along the same lines
\begin{equation}
E_k(\alpha^{\mathcal{L}(d)})=W_k^{(t)}(\alpha;d)
\underset{k\to \infty}{\sim} \frac{\Lambda^{(t)}(\alpha;d)-  \Lambda^{(t)}(\alpha;d-1)+\sqrt{8\sqrt{3}}}{\sqrt{8\sqrt{3}}}
\label{eq:avkinft}
\end{equation}
with $\Lambda^{(t)}(\alpha;d)$ fixed by
\begin{equation*}
\frac{\left(\Lambda^{(t)}(\alpha;d)\right)^2+\sqrt{8\sqrt{3}} (2 d+3) \Lambda^{(t)}(\alpha;d) +8\sqrt{3} \, d (d+3)}{\left(\Lambda^{(t)}(\alpha;d)\right)^2+\sqrt{8\sqrt{3}} (2
   d+3) \Lambda^{(t)}(\alpha;d) +8\sqrt{3} (d+1) (d+2)}=\alpha\, \frac{d\,  (d+3)}{(d+1)( d+2)}\ .
\end{equation*}
This immediately leads to  \eqref{eq:Wtintf}
 for finite $d$ and, in the limit $d\to \infty$, to \eqref{eq:CLG} with $c=1/2$.
We have similarly 
\begin{equation*}
E_k(\alpha_1^{\mathcal{L}(d_1)}\alpha_2^{\mathcal{L}(d_2)})
\underset{k\to \infty}{\sim} \frac{\Lambda^{(t)}(\alpha_1,\alpha_2;d_1,d_2)-  \Lambda^{(t)}(\alpha_1,\alpha_2;d_1-1,d_2-1)+\sqrt{8\sqrt{3}}}{\sqrt{8\sqrt{3}}}
\end{equation*}
with $\Lambda^{(t)}(\alpha_1,\alpha_2;d_1,d_2)$ now fixed by
\begin{equation*}
\begin{split}
& \hskip -1.cm \frac{\left(\Lambda^{(t)}(\alpha_1,\alpha_2;d_1,d_2)\right)^2+\sqrt{8\sqrt{3}} (2 d_2+3) \Lambda^{(t)}(\alpha_1,\alpha_2;d_1,d_2) +8\sqrt{3} \, d_2(d_2+3)}{\left(\Lambda^{(t)}(\alpha_1,\alpha_2;d_1,d_2)\right)^2+\sqrt{8\sqrt{3}} (2
   d_2+3) \Lambda^{(t)}(\alpha_1,\alpha_2;d_1,d_2) +8\sqrt{3} (d_2+1) (d_2+2)}= \\
   &\qquad \qquad \qquad  \alpha_2\, \frac{\left(\Lambda^{(t)}(\alpha_1;d_1)\right)^2+\sqrt{8\sqrt{3}}(2 d_2+3) \Lambda^{(t)}(\alpha_1;d_1) +8\sqrt{3} \, d_2 (d_2+3)}{\left(\Lambda^{(t)}(\alpha_1;d_1)\right)^2+\sqrt{8\sqrt{3}} (2
   d_2+3) \Lambda^{(t)}(\alpha_1;d_1) +8\sqrt{3} (d_2+1) (d_2+2)}\ .\\
   \end{split}
\end{equation*}
Setting $d_1=d$, $d_2=v\, d$ and sending $d\to \infty$ leads again to \eqref{eq:zuniv}, now with $c=1/2$.
This ends the proof of our first and third main results for triangulations.

\subsection{The hull perimeter at a finite fraction of the total distance}
\label{sec:finitefraction}
We end our calculation with the derivation of \eqref{eq:Ktau}, corresponding to the limit $d,k \to \infty$ with $u=d/k$ fixed ($0<u<1$).
For quadrangulations, we have
 \begin{equation*}
 \lim_{k\to \infty} E_k (e^{-\tau L(k\, u)})=  \lim_{k\to \infty} W_k^{(q)}(e^{-\frac{\tau}{(k\, u)^2}};k\, u)= \lim_{k\to \infty}  
\frac{A_3^{(q)}(e^{-\frac{\tau}{(k\, u)^2}};k\, u,k)}{A_3^{(q)}(1;k\, u,k)}
 \end{equation*}
 (recall that the denominator in the last expression is actually independent of $u$). Again the full knowledge of $A_3^{(q)}(\alpha;d,k)$ 
 is not required to handle the limit of large $k$ and we may again have recourse to the scaling limit by setting 
$k=K/\epsilon$ and letting $\epsilon\to 0$. We have now the expansions
\begin{equation}
\begin{split}
&\lambda^{(q)}\left(e^{-\frac{\tau}{(K\, u)^2} \epsilon^2};\frac{K\, u}{\epsilon}\right)=\mu^{(q)} -\nu^{(q)}\epsilon+O(\epsilon^2)\\
&\lambda^{(q)}\left(e^{-\frac{\tau}{(K\, u)^2} \epsilon^2};\frac{K\, u}{\epsilon}-1\right)=\mu^{(q)} -\xi^{(q)}\epsilon+O(\epsilon^2)\\
\end{split}
 \label{eq:explambda}
\end{equation}
where $\mu^{(q)}\equiv \mu^{(q)}(\tau;K\, u)$ is obtained by expanding \eqref{eq:lambdaddef} at second order in $\epsilon$, while  $\nu^{(q)}\equiv \nu^{(q)}(\tau;K\, u)$ and $\xi^{(q)}\equiv \xi^{(q)}(\tau;K\, u)$ follow from an expansion to third order. We find explicitly 
\begin{equation*}
\begin{split}
& \hskip -1.cm \mu^{(q)}=e^{\sqrt{6} w}\times \frac{9 w^2 \cosh ^2\left(\sqrt{\frac{3}{2}} w\right)\!+\!\left(9
   w^2\!-\!6 w \sqrt{9 w^2 \coth ^2\left(\sqrt{\frac{3}{2}} \right)\!+\!2 \tau } \!+\!2 \tau
   \right) \sinh ^2\left(\sqrt{\frac{3}{2}} w\right)}{9 w^2\!+\!2 \tau  \sinh
   ^2\left(\sqrt{\frac{3}{2}} w\right)}\\
   &
  \qquad \hbox{with}\  w=K\, u\\
  \end{split}
\end{equation*}
and
\begin{equation*}
\hskip -1.2cm \nu^{(q)}-\xi^{(q)}=\sqrt{6} \left(\frac{\left(e^{\sqrt{6} K u}+1\right) \left(e^{\sqrt{6} K u}-\mu^{(q)}
   \right)^3}{\left(e^{\sqrt{6} K u}-1\right)^3 \left(e^{\sqrt{6} K
   u}+\mu^{(q)}\right)}-\mu^{(q)} \right)\ .
   \end{equation*}
Plugging the expansion \eqref{eq:explambda} in \eqref{eq:result}, we then have the scaling behavior
\begin{equation}
\begin{split}
&Z^{(q)}\left(e^{-\frac{\tau}{(K\, u)^2} \epsilon^2};\frac{K\, u}{\epsilon},\frac{K}{\epsilon}\right)\sim \zeta^{(q)}(\tau;K,u)\,  \epsilon^3 \ , \\ 
&\zeta^{(q)}(\tau;K,u)=  \left(  \frac{24 e^{\sqrt{6} K} \left(e^{\sqrt{6} K}+\mu^{(q)}\right) \left( \nu^{(q)}- \xi^{(q)}+\sqrt{6}\, \mu^{(q)}\right)}{\left(e^{\sqrt{6} K}-\mu^{(q)}\right)^3}\right)\\
&\hskip 2.cm 
= \frac{24 \sqrt{6} \, e^{\sqrt{6} K} \left(e^{\sqrt{6} K u}+1\right) \left(e^{\sqrt{6}
   K u}-\mu^{(q)} \right)^3 \left(e^{\sqrt{6} K}+\mu^{(q)} \right)}{\left(e^{\sqrt{6} K
   u}-1\right)^3 \left(e^{\sqrt{6} K}-\mu^{(q)} \right)^3 \left(e^{\sqrt{6} K
   u}+\mu^{(q)} \right)}  \ .\\
\end{split} 
\label{eq:zeta}
\end{equation}
To be consistent with the expansion \eqref{eq:Zexpan}, we need that 
\begin{equation*}
A_i^{(q)}(e^{-\frac{\tau}{(k\, u)}};k\, u,k)\underset{k\to \infty}{\sim} k^{2i-3} \times [K^{2i-3}]\zeta^{(q)}(\tau;K,u)\ .
\end{equation*}
We deduce in particular
\begin{equation*}
  \lim_{k\to \infty} E_k (e^{-\tau L(k\, u)})=  \lim_{k\to \infty}  
\frac{A_3^{(q)}(e^{-\frac{\tau}{(k\, u)^2}};k\, u,k)}{A_3^{(q)}(1;k\, u,k)}= \frac{[K^{3}]\zeta^{(q)}(\tau;K,u)}{[K^{3}]\zeta^{(q)}(0;K,u)}\ .
\end{equation*}
Using the above explicit expressions for $\mu^{(q)}$ and $\zeta^{(q)}(\tau;K,u)$ and performing a small $K$ expansion to
extract the appropriate term of order $K^3$, we eventually arrive at \eqref{eq:Ktau} with $c=1/3$.

We can repeat the analysis for triangulations. If we use notations which parallel those introduced for quadrangulations, we find the simple 
correspondence
\begin{equation*}
\begin{split}
&\hskip -1.2cm \mu^{(t)}(\tau;K\, u)= \mu^{(q)}\left(\frac{3\tau}{2};\frac{2}{3^{1/4}}\, K\, u \right)\\
&\hskip -1.2cm \nu^{(t)}-\xi^{(t)}=\sqrt{8\sqrt{3}} \left(\frac{\left(e^{\sqrt{8\sqrt{3}} K u}+1\right) \left(e^{\sqrt{8\sqrt{3}} K u}-\mu^{(t)}
   \right)^3}{\left(e^{\sqrt{8\sqrt{3}} K u}-1\right)^3 \left(e^{\sqrt{8\sqrt{3}} K
   u}+\mu^{(t)}\right)}-\mu^{(t)} \right)\ .\\
   &\hskip -1.2cm  \zeta^{(t)}(\tau;K,u)= \frac{24\, \sqrt{8\sqrt{3}} \, e^{\sqrt{8\sqrt{3}} K} \left(e^{\sqrt{8\sqrt{3}} K u}+1\right) \left(e^{\sqrt{8\sqrt{3}}
   K u}-\mu^{(t)} \right)^3 \left(e^{\sqrt{8\sqrt{3}} K}+\mu^{(t)} \right)}{\left(e^{\sqrt{8\sqrt{3}} K
   u}-1\right)^3 \left(e^{\sqrt{8\sqrt{3}} K}-\mu^{t)} \right)^3 \left(e^{ \sqrt{8\sqrt{3}} K
   u}+\mu^{(t)} \right)} \\
\end{split}
\end{equation*}
and
\begin{equation*} 
\lim_{k\to \infty} E_k (e^{-\tau L(k\, u)})= \frac{[K^{3}]\zeta^{(t)}(\tau;K,u)}{[K^{3}]\zeta^{(t)}(0;K,u)}\ .\\
\end{equation*}
With these expressions, we arrive at the same expression \eqref{eq:Ktau} as for quadrangulations, but now with $c=1/2$\footnote{Note that the
passage from quadrangulations to triangulations amounts to two simultaneous changes: $K\to 2 K/3^{1/4}$ and $\tau\to 3\tau /2$. The first
change has no impact on the formula for $\lim_{k\to \infty} E_k (e^{-\tau L(k\, u)})$ but the second change is responsible for
the passage from $c=1/3$ to $c=1/2$.} .

This proves our second main result for quadrangulations and triangulations. To obtain the expression \eqref{eq:Pk} for $\mathcal{P}(L;u)$, we simply have to compute the inverse Laplace transform 
of $F\left(\sigma(\tau;u);u\right)$ as given by \eqref{eq:Ktau}. The details of this computation
are presented in Appendix A.  

\section{Conclusion}
\label{sec:conclusion}
In this paper, we showed how to control the hull perimeters in pointed-rooted quadrangulations or triangulations 
in a very explicit way by some appropriate decoration of recursion relations inherited from 
a decomposition of the maps via cuts appearing precisely 
along hull boundaries. Even though we concentrated here only 
on the statistics of hull perimeters, namely the lengths of hull boundaries, we could in principle measure other quantities characterizing
the hulls such as for instance their volume, as was done in \cite{CLG14a,CLG14b}. 

It was recognized in \cite{Krikun03,Krikun05} that the structure of the hull in a random map can be understood as some particular 
time-reversed branching process. In our formalism, this process appears in the branching nature of the successive cuttings
of the original slice encoding the pointed-rooted map at hand into smaller and smaller sub-slices. In particular, the branching information 
is entirely captured by the kernel of the recursion relation. In this respect, it is likely that the universal laws that we found can be 
given some more direct interpretation as statistical laws for appropriate quantities in the branching process.  

To conclude, we would like to stress that we were eventually interested in this paper in the limit $N\to\infty$ with distances 
which do not scale with $N$.
This is the so-called local limit of large maps, whose continuous description is provided by the Brownian plane \cite{CLG13}.
If we want
instead a full access to properties of the so-called scaling limit, described by the celebrated Brownian map \cite{Miermont2013,legall2013}, we
simply have to let $k$ and $d$ scale like $N^{1/4}$ when $N$ becomes large.  Our expressions are also well adapted
to this scaling limit and it could be interesting to extend our calculations to this case.
 
Note finally that the notion of hull was recently shown in \cite{MS15} to be a fundamental ingredient in the characterization
of the Brownian map, which makes its statistical study even more appealing.

\appendix
\section{Inverse Laplace transforms}
The quantity $\mathcal{P}(L;u)$ is simply the inverse Laplace transform 
of $F\left(\sigma(\tau;u);u\right)$ as given by \eqref{eq:Ktau}, where $L$ is the conjugate variable of $\tau$.
Assuming that we know the inverse Laplace transform  $\check{F}(\ell;u)$ of $F(\sigma;u)$, where $\ell$ in the conjugate variable of $\sigma$, we clearly 
have, since $\sigma(\tau;u)=(b(u)-1)+c\, b(u)\times \tau$ (with $b(u)=(1-u)^2/u^2$):
\begin{equation*}
\mathcal{P}(L;u)= \frac{1}{c\, b(u)}\,  \check{F}\left(\frac{L}{c\, b(u)};u\right)\, e^{-\frac{b(u)-1}{c\, b(u)}\, L}\ .
\end{equation*}
To compute $\check{F}(\ell;u)$ , we simply need the inverse Laplace transforms of the functions $\sigma^n$ for $n=-4$, $-3$ and $-2$,
and those of $\sigma^n/(1+\sigma)^{5/2}$ for $n=-4,-3,\cdots, +1$. The first three are respectively $\ell^3/3!$, $\ell^2/2!$ and $\ell$.
As for the last six, they are equal to:
\begin{equation*}
\begin{split}
&n=-4:\ \ \qquad  \frac{e^{-\ell} \sqrt{\ell} \left(4
   \ell^2+315\right)}{24 \sqrt{\pi }}+\frac{1}{48} \left(8 \ell^3-60 \ell^2+210 \ell-315\right)\,
   \text{erf}\left(\sqrt{\ell}\right)\ ,\\& n=-3:\qquad  -\frac{5 e^{-\ell} \sqrt{\ell} (2 \ell+21)}{12 \sqrt{\pi
   }}
   +\frac{1}{8} (4\ell^2-20 \ell+35)\, \text{erf}\left(\sqrt{\ell}\right)\ ,\\
   & n=-2:\ \ \qquad   \frac{e^{-\ell} \sqrt{\ell} (4 \ell+15)}{3 \sqrt{\pi }}+\frac{1}{2}\left(2\ell-5\right)\, \text{erf}\left(\sqrt{\ell}\right)\ ,\\ 
   & n=-1 :\qquad -\frac{2 e^{-\ell}
   \sqrt{\ell} (2 \ell+3)}{3 \sqrt{\pi }}+\text{erf}\left(\sqrt{\ell}\right)\ ,\\& n=0 :\qquad  \ \ \ \ \frac{4 e^{-\ell} \ell^{3/2}}{3 \sqrt{\pi
   }}\ ,\\&n=1 :\qquad \ \ - \frac{2 e^{-\ell}  \sqrt{\ell} (2 \ell-3) }{3 \sqrt{\pi }}\ .\\
\end{split}
\end{equation*}
Using these values, we easily evaluate 
\begin{equation*}
\check{F}(\ell;u)=\frac{u^3 \sqrt{\ell}\, e^{-\ell} }{2
   \sqrt{\pi }} \left(\left(e^\ell \sqrt{\pi \, \ell}
   \left(1-\text{erf}\left(\sqrt{\ell}\right)\right)-1\right) p(\ell)+r(\ell)\right)
\end{equation*} 
with $p(\ell)$ and $r(\ell)$ as in \eqref{eq:Pk} where $b=b(u)=(1-u)^2/u^2$. Equation \eqref{eq:Pk} follows immediately.
\vskip .5cm
Let us know discuss how we may obtain the expression \eqref{eq:PL1L2} for the joint probability
density $\mathcal{P}(L_1,L_2;v)$. Clearly, it is simply the result of a double inverse Laplace transform
of the expression \eqref{eq:zuniv} where $L_1$ and $L_2$ are the conjugate variables of $\tau_1$ and $\tau_2$ respectively.
The inverse Laplace transform in the variable $L_2$ is easily performed, leading to
\begin{equation*}
\begin{split}
& 2   \sqrt{\frac{L_2}{\pi}}\, \frac{v^3}{c^{3/2}}\frac{e^{-\frac{L_2\, v^2 \omega _1^2}{c \left((v-1) \omega _1+1\right)^2}}
  }{  \left((v-1) \omega_1+1\right)^3}
 =\sqrt{\frac{2}{\pi }}\, \frac{(v-1) v^2}{c}\, \times 
 \frac{ e^{-\frac{t_2^2 \, (v-1)^2  \omega_1^2}{2 \left((v-1) \omega_1+1\right)^2}} \, t_2}{\left((v-1)
   \omega_1+1\right)^3}\\
   & \hbox{with}\quad  \omega_1\equiv \sqrt{1+c\,  \tau_1}\quad  \hbox{and}\quad   t_2\equiv  \sqrt{\frac{2 L_2\, v^2}{c\, (v-1)^2}}\ . \\
   \end{split}
\end{equation*}   
Writing the second factor as 
\begin{equation*}
    \frac{ e^{-\frac{t_2^2 \, (v-1)^2  \omega_1^2}{2 \left((v-1) \omega_1+1\right)^2}} \, t_2}{\left((v-1)
   \omega _1+1\right)^3}= \frac{e^{-\frac{t_2^2}{2 \left(\Omega_1+1\right)^2}} t_2 \Omega_1^3}{\left(\Omega_1+1\right)^3}\ , \qquad \Omega_1\equiv \frac{1}{(v-1) \omega_1}\ , 
   \end{equation*}
 we may expand the last expression as
 \begin{equation*}
 \frac{e^{-\frac{t_2^2}{2 \left(\Omega_1+1\right)^2}} t_2 \Omega _1^3}{\left(\Omega_1+1\right)^3}=e^{-\frac{t_2^2}{2}}\times \sum_{n=0}^{\infty} \frac{\Omega_1^{n+3}}{n!} (-1)^n \, \pi_n(t_2)
 \end{equation*}
 with $\pi_n(t)$ as in \eqref{eq:pin}. This latter equation 
 is obtained from the relation\footnote{Using the Leibniz formula to compute the $n$-th derivative and rearranging the terms, the
 right hand side in the relation is indeed easily rewritten as the shift operator $\sum_{n=0}^\infty \frac{a^n}{n!}\frac{\partial^n}{\partial t^n}f(t)=f(t+a)$
 acting on the function $f(t)=e^{-t^2/2}$ with $a=-t \, \Omega_1/(\Omega_1+1)$.}
 \begin{equation*}
 e^{-\frac{t_2^2}{2(\Omega_1+1)^2}}=\sum_{n=0}^\infty \frac{(-\Omega_1)^n}{n!} t\, \frac{\partial}{\partial t^n}\left(t^{n-1}e^{-\frac{t_2^2}{2}}\right)
 \end{equation*}
 upon differentiating with expect to $\Omega_1$ and multiplying by $(\Omega_1^3/t_2)$.
 We may now use the inverse Laplace transform:
 \begin{equation*}
 \hskip -1.cm \Omega_1^{n+3}=\frac{1}{(v-1)^{n+3}}\frac{1}{(1+c\, \tau_1)^{\frac{n+3}{2}}}\quad \xrightarrow[]{\text{Inv. Laplace Trans.}}\quad 
 \frac{e^{-\frac{L_1}{c}}    \left(\frac{L_1}{c (v-1)^2}\right)^{\frac{n+3}{2}-1}}{c\, (v-1)^2 \, \Gamma
   \left(\frac{n+3}{2}\right)}\ .
 \end{equation*}
 Incorporating this formula in our infinite sum, we end up with
 \begin{equation*}
\hskip -1.cm \mathcal{P}(L_1,L_2;v)= \sqrt{\frac{2}{\pi }}\, \frac{(v-1) v^2}{c}\times  e^{-\frac{t_2^2}{2}} \times \sum_{n=0}^{\infty} \frac{1}{n!} 
  \frac{e^{-\frac{L_1}{c}}    \left(\frac{L_1}{c (v-1)^2}\right)^{\frac{n+3}{2}-1}}{c\, (v-1)^2 \, \Gamma
   \left(\frac{n+3}{2}\right)}
 (-1)^n \, \pi_n(t_2)\ ,
 \end{equation*}
 which is precisely \eqref{eq:PL1L2}.
 
 \section{First terms of the $g$-expansion of $Z^{(q)}(\alpha;d,k)$ and $Z^{(t)}(\alpha;d,k)$}
We give here the first terms of the $g$-expansion of $Z^{(q)}(\alpha;d,k)$ (up to order $g^8$) and of 
$Z^{(t)}(\alpha;d,k)$ (up to order $g^{12}$) for the smallest values of $d$ and $k$ allowed in our formulas. 
The coefficients of $\alpha^{\mathcal{L}}\, g^N$ in these expansions are the \emph{numbers} of $k$-pointed-rooted quadrangulations 
(respectively triangulations) with a given number $N$ of faces and a fixed value $\mathcal{L}$ of the hull 
perimeter $\mathcal{L}(d)$. We find:
\begin{equation*}
\begin{split}
& \hskip -1.cm Z^{(q)}(\alpha;2,3)= 
\alpha ^2 g^2\!+\!15 \alpha ^2 g^3\!+\!\alpha ^2 \left(\alpha ^2\!+\!178\right) g^4\!+\!7 \alpha
   ^2 \left(4 \alpha ^2\!+\!281\right) g^5\\
   & \hskip 2.cm \!+\!\alpha ^2 \left(\alpha ^4\!+\!518 \alpha
   ^2\!+\!21165\right) g^6\!+\!\left(42 \alpha ^6\!+\!8018 \alpha ^4\!+\!225488 \alpha ^2\right)
   g^7\\
   & \hskip 2.cm \!+\!\alpha ^2 \left(\alpha ^6\!+\!1075 \alpha ^4\!+\!112671 \alpha ^2\!+\!2395983\right)
   g^8\!+\!O\left(g^9\right)\\
& \hskip -1.cm Z^{(q)}(\alpha;2,4)=
\alpha ^2 g^3\!+\!22 \alpha ^2 g^4\!+\!\left(\alpha ^4\!+\!342 \alpha ^2\right) g^5\!+\!\left(36
   \alpha ^4\!+\!4640 \alpha ^2\right) g^6\\
   &\hskip 2.cm \!+\!\left(\alpha ^6\!+\!815 \alpha ^4\!+\!58799 \alpha
   ^2\right) g^7\!+\!\left(51 \alpha ^6\!+\!14914 \alpha ^4\!+\!716865 \alpha ^2\right)
 g^8\!+\!O\left(g^9\right)
\\
& \hskip -1.cm Z^{(q)}(\alpha;3,4)=
\alpha ^2 g^3\!+\!22 \alpha ^2 g^4\!+\!\alpha ^2 \left(2 \alpha ^2\!+\!341\right) g^5\!+\!\alpha
   ^2 \left(71 \alpha ^2\!+\!4605\right) g^6\\
   &\hskip 2.cm \!+\!\alpha ^2 \left(3 \alpha ^4\!+\!1586 \alpha
   ^2\!+\!58026\right) g^7\!+\!5 \alpha ^2 \left(30 \alpha ^4\!+\!5731 \alpha ^2\!+\!140605\right)
   g^8\!+\!O\left(g^9\right)
\\
& \hskip -1.cm Z^{(q)}(\alpha;2,5)=
\alpha ^2 g^4\!+\!29 \alpha ^2 g^5\!+\!\left(\alpha ^4\!+\!555 \alpha ^2\right) g^6\!+\!\left(43
   \alpha ^4\!+\!8876 \alpha ^2\right) g^7\\
   &\hskip 2.cm \!+\!\left(\alpha ^6\!+\!1127 \alpha ^4\!+\!128712
   \alpha ^2\right) g^8\!+\!O\left(g^9\right)
\\
& \hskip -1.cm Z^{(q)}(\alpha;3,5)=
\alpha ^2 g^4\!+\!29 \alpha ^2 g^5\!+\!2 \alpha ^2 \left(\alpha ^2\!+\!277\right) g^6\!+\!\left(87
   \alpha ^4\!+\!8832 \alpha ^2\right) g^7\\
   &\hskip 2.cm \!+\!\alpha ^2 \left(3 \alpha ^4\!+\!2300 \alpha
   ^2\!+\!127537\right) g^8\!+\!O\left(g^9\right)
\\
& \hskip -1.cm Z^{(q)}(\alpha;4,5)=
\alpha ^2 g^4\!+\!29 \alpha ^2 g^5\!+\!2 \alpha ^2 \left(\alpha ^2\!+\!277\right) g^6\!+\!\alpha
   ^2 \left(86 \alpha ^2\!+\!8833\right) g^7\\
   &\hskip 2.cm \!+\!\alpha ^2 \left(3 \alpha ^4\!+\!2251 \alpha
   ^2\!+\!127586\right) g^8\!+\!O\left(g^9\right)
\\
\end{split}
\end{equation*}
and
\begin{equation*}
\begin{split}
& \hskip -1.cm Z^{(t)}(\alpha;1,2)= 
\alpha  g^2\!+\!\alpha  (\alpha \!+\!14) g^4\!+\!\alpha  \left(\alpha ^2\!+\!26 \alpha \!+\!199\right) g^6\!+\!\alpha  \left(\alpha ^3\!+\!39 \alpha ^2\!+\!533 \alpha
   \!+\!2952\right) g^8 \\ &
   \hskip 3.cm \!+\!\alpha  \left(\alpha ^4\!+\!53 \alpha ^3\!+\!1062 \alpha ^2\!+\!10147 \alpha \!+\!45473\right) g^{10}
   \\ & \hskip 3.cm \!+\!\alpha  \left(\alpha ^5\!+\!68
   \alpha ^4\!+\!1824 \alpha ^3\!+\!25040 \alpha ^2\!+\!187756 \alpha \!+\!722498\right) g^{12}\!+\!O\left(g^{14}\right)
\\
& \hskip -1.cm Z^{(t)}(\alpha;1,3)=
\alpha  g^4\!+\!\alpha  (\alpha \!+\!28) g^6\!+\!\alpha  \left(\alpha ^2\!+\!42 \alpha \!+\!612\right) g^8\!+\!\alpha  \left(\alpha ^3\!+\!57 \alpha ^2\!+\!1220
   \alpha \!+\!12326\right) g^{10}\\
   & \hskip 4.cm \!+\!\alpha  \left(\alpha ^4\!+\!73 \alpha ^3\!+\!2090 \alpha ^2\!+\!30456 \alpha \!+\!239793\right) g^{12}\!+\!\alpha 
\!+\!O\left(g^{14}\right)
\\
& \hskip -1.cm Z^{(t)}(\alpha;2,3)=
\alpha  g^4\!+\!\alpha  (2 \alpha \!+\!27) g^6\!+\!\alpha  \left(3 \alpha ^2\!+\!79 \alpha \!+\!573\right) g^8\!+\!\alpha  \left(4 \alpha ^3\!+\!159 \alpha
   ^2\!+\!2178 \alpha \!+\!11263\right) g^{10}\\ 
   & \hskip 4.cm \!+\!\alpha  \left(5 \alpha ^4\!+\!270 \alpha ^3\!+\!5479 \alpha ^2\!+\!51970 \alpha \!+\!214689\right)
   g^{12}\!+\!O\left(g^{14}\right)
\\
& \hskip -1.cm Z^{(t)}(\alpha;1,4)=
\alpha  g^6\!+\!\alpha  (\alpha \!+\!42) g^8\!+\!\alpha  \left(\alpha ^2\!+\!56 \alpha \!+\!1225\right) g^{10}\\
& \hskip 4.cm \!+\!\alpha  \left(\alpha ^3\!+\!71 \alpha ^2\!+\!2031
   \alpha \!+\!30792\right) g^{12}\!+\!O\left(g^{14}\right)
\\
& \hskip -1.cm Z^{(t)}(\alpha;2,4)=
\alpha  g^6\!+\!\alpha  (2 \alpha \!+\!41) g^8\!+\!\alpha  \left(3 \alpha ^2\!+\!111 \alpha \!+\!1168\right) g^{10}\\
& \hskip 4.cm \!+\!\alpha  \left(4 \alpha ^3\!+\!213 \alpha
   ^2\!+\!3984 \alpha \!+\!28694\right) g^{12}\!+\!O\left(g^{14}\right)
\\
& \hskip -1.cm Z^{(t)}(\alpha;3,4)=
\alpha  g^6\!+\!\alpha  (2 \alpha \!+\!41) g^8\!+\!\alpha  \left(3 \alpha ^2\!+\!108 \alpha \!+\!1171\right) g^{10}\\
& \hskip 4.cm \!+\!\alpha  \left(4 \alpha ^3\!+\!204 \alpha
   ^2\!+\!3791 \alpha \!+\!28896\right) g^{12}\!+\!O\left(g^{14}\right)\ .
\\
\end{split}
\end{equation*}

\section{General expression for the average hull perimeter}
\begin{figure}
\begin{center}
\includegraphics[width=8.5cm]{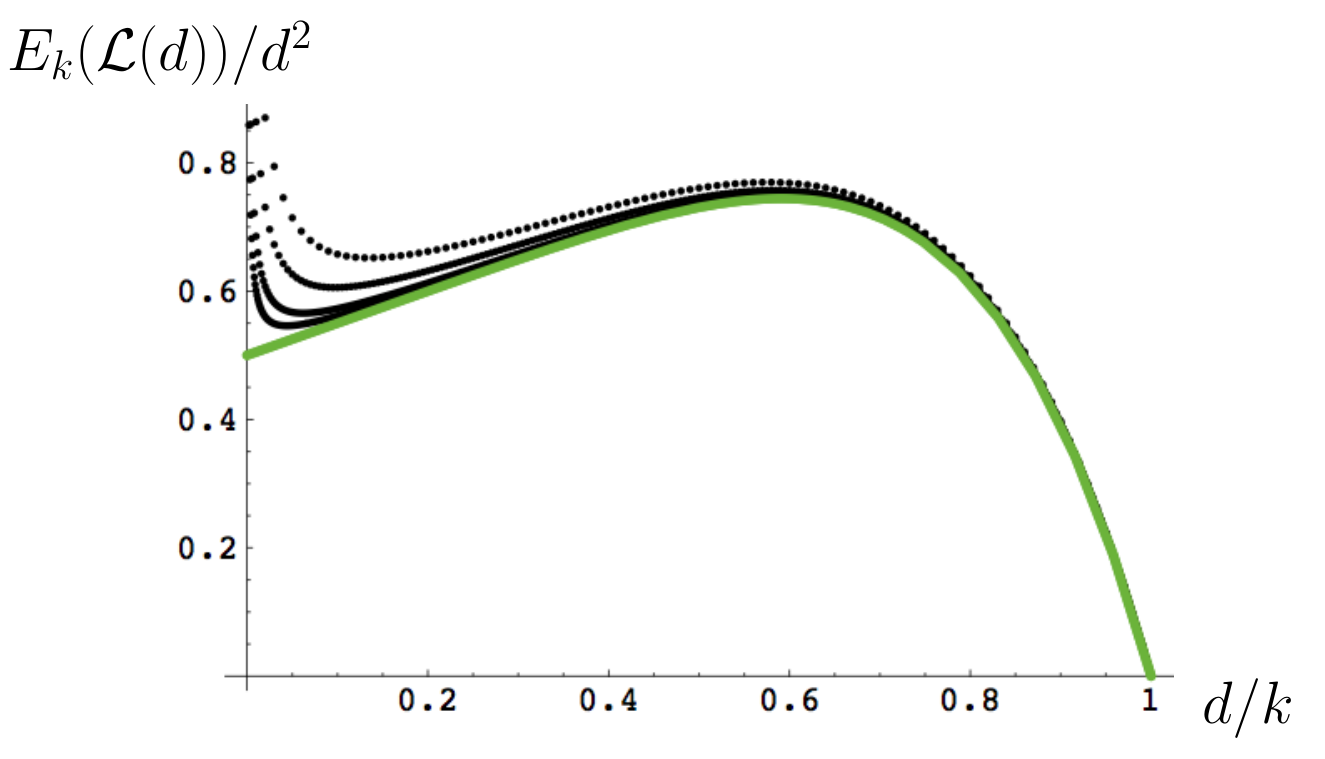}
\end{center}
\caption{A plot of $E_k(\mathcal{L}(d))/d^2$ for quadrangulations as a function of $d/k$ for $k=100$, $200$, $500$ and
$1000$ (from top to bottom), and the scaling limit $L_{\rm av}(d/k)$ (green line).}
\label{fig:finiteLkd}
\end{figure}

We give here for information the expression of $E_k(\mathcal{L}(d))$ at
\emph{finite $k$ and $d$}, 
for both quadrangulations (satisfying $E_k(\mathcal{L}(d))=\frac{\partial}{\partial \alpha} W_k^{(q)}(\alpha;d)\vert_{\alpha=1}$)
and triangulations ($E_k(\mathcal{L}(d))=\frac{\partial}{\partial \alpha} W_k^{(t)}(\alpha;d)\vert_{\alpha=1}$).
In the case of quadrangulations, we find:
\begin{equation*}
\begin{split}
& \hskip -1.2cm E_k(\mathcal{L}(d))=
\frac{ k (k\!+\!1) (k\!+\!2)}{2 \left(k^2\!+\!2 k\!-\!1\right) \left(5 k^4\!+\!20 k^3\!+\!27 k^2\!+\!14
   k\!+\!4\right)}
   \times \\
&  \hskip -1.2cm\times \Bigg(
(d\!-\!1) (d\!+\!1) (d\!+\!2) (d\!+\!4) (2 k\!+\!3)\times \\ & \times \frac{ \left((1\!-\!d) (d\!+\!1) (d\!+\!2) (d\!+\!4) \left(5 d^2\!+\!15 d\!+\!17\right)\!+\!(k\!+\!1)^2 (k\!+\!2)^2 \left(5 k^2\!+\!15 k\!+\!2\right)\!-\!4\right)}{3 (2 d\!+\!3)
   (k\!+\!1)^2 (k\!+\!2)^2}
\\
   & \hskip -.8cm-
(d\!-\!2) d (d\!+\!1) (d\!+\!3) (2 k\!+\!1) \times \\  & \times \frac{\left((2\!-\!d) d (d\!+\!1) (d\!+\!3) \left(5 d^2\!+\!5 d\!+\!7\right)\!+\!k^2 (k\!+\!1)^2 \left(5 k^2\!+\!5 k\!-\!8\right)\!-\!4\right)}{3(2 d\!+\!1) k^2 (k\!+\!1)^2}  
   \Bigg)\ .\\
   \end{split}
\end{equation*}
for $2\leq d<k$. This quantity, rescaled by $d^2$, is plotted in figure \ref{fig:finiteLkd} to emphasize its scaling limit at large $k$ and $d$.

In the case of triangulations, we find instead:
\begin{equation*}
\begin{split}
& \hskip -1.2cm E_k(\mathcal{L}(d))=
\frac{k^2 (k\!+\!1)^2 }{2(2 k\!+\!1) \left(5 k^6\!+\!15 k^5\!+\!14 k^4\!+\!3 k^3\!-\!k^2\!-\!1\right)}\times \\
&  \hskip -1.2cm\times \left(\frac{d (d\!+\!1) (d\!+\!2) (d\!+\!3) \left(10 (k\!+\!1)^6\!-\!7(k\!+\!1)^4\!-\!2 d (d\!+\!1) (d\!+\!2) (d\!+\!3) \left(5 d^2\!+\!15 d\!+\!14\right)\!-\!2\right)}{(k\!+\!1)^3 (2 d\!+\!3)}\right .\\
   & \hskip -.8cm\left.-\frac{ (d\!-\!1) d (d\!+\!1) (d\!+\!2) \left(10 k^6\!-\!7 k^4\!-\!2 (d\!-\!1) d (d\!+\!1) (d\!+\!2) \left(5 d^2\!+\!5 d\!+\!4\right)\!-\!2\right)}{ k^3 (2 d\!+\!1)}\right)\ .\\
   \end{split}
\end{equation*}
for $1\leq d<k$.
 
\section*{Acknowledgements} 
The author acknowledges the support of the grant ANR-14-CE25-0014 (ANR GRAAL).
 
\bibliographystyle{plain}
\bibliography{hull}

\end{document}